\documentclass{article}
\usepackage{amsfonts,amsmath,theorem}
\usepackage{color}

\newtheorem{theorem}{Theorem}[section]

\newtheorem{remark}{Remark}[section]
\newtheorem{prop}{Proposition}[section]
\newtheorem{proposition}{Proposition}[section]
\newtheorem{defi}{Definition}[section]
\newcommand{\bu}{$\bullet$}

\newcommand{\nn}{\nonumber}

\newcommand{\be}{\begin{equation}}
\newcommand{\ee}{\end{equation}}
\newcommand{\bea}{\begin{eqnarray}}
\newcommand{\eea}{\end{eqnarray}}
\newcommand{\bd}{\begin{displaymath}}
\newcommand{\ed}{\end{displaymath}}

\def \R{{\mathbf R}}

\def \Re{Re}                    

\def \R {\mathbf R}
\def \Z {\mathbf Z}

\def \eps{{\epsilon}}

\def \uT {\tilde u}

\def \del {{\partial}}

\def\sgn {\hbox{sign}}

{\title{Mathematics and Turbulence: where do we stand?}}

\author{
Claude BARDOS%
\footnote{
Universit\'e Denis Diderot  and Laboratoire J.-L. Lions, Universit\'e Pierre et Marie
Curie, Paris, France (claude.bardos@gmail.com).},
and
Edriss S. TITI%
\footnote{
Department of Computer Science and Applied Mathematics,
Weizmann Institute of Science,
Rehovot 76100, Israel;
{\bf also}
Department of Mathematics
and  Department of Mechanical and  Aerospace Engineering,
University of California
Irvine, CA  92697-3875, USA;
(etiti@math.uci.edu). }}%

\date{January 29, 2013}
\begin{document}
\maketitle
\begin{abstract}
This contribution covers the topics presented by the   authors at the  {\it ``Fundamental Problems of Turbulence, 50 Years after the Marseille Conference 1961"} meeting that took place in Marseille in 2011.
It  focuses on some of the mathematical approaches to fluid dynamics and turbulence. This  contribution does not pretend to cover or answer, as the reader may discover,  the fundamental  questions  in turbulence, however, it aims toward presenting some of the most recent advances in attacking these questions using rigorous mathematical tools. Moreover, we consider that the proofs  of the mathematical statements  (concerning, for instance,  finite time regularity, weak solutions and vanishing viscosity) may contain  information  as relevant, to the understanding of the underlying problem, as    the statements themselves.
\end{abstract}
\maketitle
\noindent {\bf MSC Subject Classifications:} 35Q30, 35Q31, 76D03, 76D09.

\noindent {\bf Keywords:} Euler equations, Navier-Stokes equations, turbulence, viscosity limit.

\section{Introduction}

One of the main goals of this contribution is to underline the relationships between, on the one hand, standard mathematical questions and, on the other hand,  classical issues in turbulence theory.  However, it is important to emphasize that there are many different mathematical approaches aiming toward understanding the complex phenomenon of turbulence,   and that only some of these aspects are stressed in this contribution (see also the other contributions in this special issue; and for  a general  user friendly introduction (for mathematicians and engineers)  the book of Frisch \cite{FR}).

A nice survey of the dynamical systems approach to turbulence and statistical solutions of Navier-Stokes equations is presented in the monographs \cite{CFT}, \cite{FMRT}, and in   Foias \cite{Foias}.  Moreover, rigorous conditions for the main tenets of the Batchelor-Kraichnan-Leith theory of 2D turbulence, namely the energy power law, enstrophy cascade, and dissipation law
can be found in \cite{BFJR}, \cite{DFJ}, \cite{FJMR} and references therein. Furthermore, recent interesting surveys concerning analytical studies of Euler equations can be found in \cite{Bardos-Titi-survey}, \cite{Constantin}, \cite{ConstantinAMS}, \cite{Gibbon} and references therein (see also  \cite{BustKerr}, \cite{Gibbon},\cite{Hou}  and references therein for  recent surveys  concerning computational studies of Euler equations).

Starting from the macroscopic Navier-Stokes-Fourier system for compressible fluid flows, we will recall the derivation of the equations of motion for incompressible fluid flows, by taking the small Mach number limit. Then by another  new scaling argument, which relies on the introduction of the Reynolds number, we will introduce the Euler equations as the formal zero viscosity limit of the Navier-Stokes equations. Historically, the Euler equations for ideal fluids were the first to be derived (Euler 1757), viscous effects being added only later to give the  Navier-Stokes equations. These viscous effects are, as it will be stressed in the paper, essential for  understanding  the interaction between solid bodies and fluid flows. Moreover, recent mathematical  results indicate that a more profound understanding of  Euler equations   seems to be essential for making any significant  progress in this field. Therefore, section \ref{euler0} is devoted to the mathematical theory of  Euler equations with, on the one hand, emphasis on the classical local well-posedness of the problem in the class of smooth solutions and, on the other hand, on the more recent notions of weak solutions. We will then return, {{in section 4}}, to the Navier-Stokes equations, mostly to consider the boundary value problems and the {{behavior}} of the solutions, in particular boundary effects, when the viscosity tends to zero.

\section{Preliminaries}\label{preliminaries}

Our approach follows three  threads. First, the notion of the {\it Cauchy problem} for an underlying evolution equation, which means the study of the corresponding solutions of this equation, for $t>t_0$, given their state  at time $t=t_0$ (often $t_0$ is taken equal to $0$).  Here the relevant issues are:  the stability or instability of the solutions, their qualitative long-time behavior, and eventually  the reliability of this equation for modeling the underlying phenomena at the relevant scales   and under the appropriate conditions.

The second thread is the notion of {\it weak solutions}. This idea preexisted in mechanics and physics in the context of, for instance,  calculus of variations; however,  in the previous century it became systematically used with the development of functional analysis and the introduction, by Gelfand, Schwartz and Sobolev,  of  the notion of distributions (which are also called, by some people,  generalized functions).
In a nutshell, consider for example,   the finite linear system of   ordinary differential equations:
\begin{equation}\label{ODE}
\frac{d}{dt} u + A(t)u=f\,,\quad  u(0)=u_0,
\end{equation}
where $u$ is an unknown vector of $n$components and $A(t)$ is  a given $n\times n$  matrix with bounded entries. In case $A(t)$ is not continuous one does not expect $u$ to be differentiable in the classical sense. Nonetheless, one would still like to make sense of the above initial value problem. To this end,  one takes the Euclidean inner product, that we denote by $\langle \cdot, \cdot\rangle$, of  equation (\ref{ODE}) with any  vector field $\phi(t)$, member of a  family of conveniently smooth  ``test functions" that are identically  equal to zero when $t$ becomes large. Then formally (because no regularity on $u$ is assumed)  integrating by parts one can consider instead the equation

\begin{equation}\label{ODE-weak}
\int_0^\infty\left [ - \langle u(t) , \frac{d}{dt}\phi(t) \rangle + \langle u(t), A^*(t)\phi(t)\rangle \right ]\,dt= \langle u(0) ,\phi(0) \rangle\,.
\end{equation}

A function $u(t)$ is said to be a weak solution of (\ref{ODE}) if $u$ satisfies equation (\ref{ODE-weak}) for every test function $\phi \in C^\infty([0,\infty))$, and $\phi(t)$ is identically zero for $t$ large enough. Notice that in order for $u$ to satisfy (\ref{ODE-weak}) it does not need to be a differentiable function.

In fact, without the refinement of functional analysis, this idea already goes back to Euler and Lagrange in the context of deriving the Euler-Lagrange equation in the calculus of variations. Similarly,  in the context of gas dynamics there is a good reason to introduce such a notion because  compressible fluids support,  and {\it  even when starting from smooth initial data} may develop, shock waves which are discontinuous, and therefore have no derivative in the classical sense. Since Leray did not know much about the regularity of the solutions of the Navier-Stokes equations it   was crucial in his seminal work \cite{Leray} to introduce weak solutions.   Of course, this weak formulation leads to  a loss of information which, for instance, may prevent the analysis of stability, or the proof of uniqueness.

Weak solutions are defined as elements of some  spaces of functions,  which will be used below with no further introduction.

With these prerequisites,  a  final  and third thread will be the consideration of {\it the equations of macroscopic perfect fluids and the relationship between them.}

We will start with the evolution of density, velocity and temperature $(\rho, u, \theta)$  given by the compressible Navier-Stokes-Fourier system of equations (cf. \cite{Landau-Lifshitz} and \cite{PLLions}):
\begin{equation*}
\begin{aligned}
& \del_t \rho  + \nabla \cdot (\rho   u ) =0 \,, \\
&  \rho \big( \del_t + u  \cdot \nabla \big) u
  + \nabla p
  = \epsilon \nabla \cdot [\nu  \sigma(u )]\,,\nn\\
&  \frac32 \rho  \big( \del_t + u \cdot \nabla \big) \theta
  + \rho \theta  \nabla \cdot u
  =  \frac12\epsilon  \nu  (\sigma(u ):  \sigma(u  ))
  + \epsilon \nabla [\kappa  \nabla\theta ] \,, \\
  \end{aligned}
  \end{equation*}
  with the state equations for the pressure, $p  = \rho  \theta$, and the stress tensor:
  \begin{equation*}
   \sigma(u)=\left(\nabla_xu+(\nabla _xu)^T\right)-\frac23(\nabla_x \cdot u) I\,.
 \end{equation*}
 In these equations,
 for any vector $u\in \R^d$  and any (scalar or vector valued function) $v$, the term  $(u\cdot \nabla ) v $ denotes the expression
 $$
 (u\cdot \nabla) v= \sum_{1\le j \le d} u_j\frac{\del v}{\del x_j}  \,,
 $$
while $\nu$ and $\kappa$ denote the kinematic viscosity  and the thermal diffusivity of the gas, respectively. These quantities depend on the nature of the molecular interaction, while $\epsilon$ denotes the dimensionless Knudsen number (proportional, up to a fudge constant, to the mean free path) which depends on the degree of rarefaction of the gas. A simplified system is then obtained (see, e.g., section 1.1 in \cite{PLLions}) for the  fluctuations of  order $\epsilon$  of density,  velocity and  temperature  about a constant background state, $(\overline{\rho},0,\overline{\theta})$, namely,
 \begin{equation}
 \rho= \overline{\rho} +\epsilon \tilde \rho, \quad \theta=   \overline{\theta}+ \epsilon \tilde \theta, \quad u=\epsilon \tilde u\,.
 \end{equation}
In order to observe the cumulative effect of the small velocity, one has to consider  the system for a time scale long enough (of the order of $\epsilon^{-1}$), and hence by rescaling the time, by factor of $\epsilon^{-1}$, and following, for instance, the analysis presented at section 1.1 in \cite{PLLions},  one obtains  the following system:

  \begin{eqnarray*}
 && \nabla \cdot \tilde u =0\,, \quad \nabla  (\overline{\theta}\tilde \rho+ \overline{\rho}\tilde \theta)=0\,,
 \\
&&\overline{\rho}(\del_t \tilde u +(\tilde u \cdot \nabla) \tilde u)  + \nabla \pi=\nu \Delta \tilde u\,,\\
&&\frac32\overline{\rho}(\del_t\tilde \theta+(\tilde u \cdot\nabla) \tilde \theta )=\kappa\Delta \tilde \theta \,,
\end{eqnarray*}
where $\pi$  is the unknown pressure.
This system is often referred to as the incompressible Navier-Stokes-Fourier system of equations. Observe that the equations  for the velocity are independent of the temperature, while the  advection term, $(\tilde u \cdot\nabla) \tilde \theta$, in the   temperature equation depends on the velocity.  The  term incompressible comes  from the $\nabla \cdot \tilde{u} =0$ condition, which implies that  the Lagrangian flow map defined by:
\begin{equation*}
x\mapsto \phi(x,t) \,, \quad \phi(x,0)=x\,,\quad \dot \phi(x,t)= \tilde u(\phi(x,t),t) \,,
\end{equation*}
is volume preserving.   The equation $\nabla(\overline{\theta}\tilde{\rho}  +\overline{\rho}\tilde{\theta})=0$ is known as  the Bernoulli relation. Together with the equation for the temperature they describe how the flow, at this    approximation level, propagate the  initial fluctuations of  density and  temperature.
In particular, if $\tilde \rho=\tilde \theta=0$ the system reduces (hereafter the  tilde  will be omitted, $\overline{\rho}$ will be taken equal to $1$, and following the tradition in this field we will replace the notation for the pressure $\pi$ by $p$) to the standard incompressible Navier-Stokes equations:

\begin{equation}
\del_t u + (u\cdot \nabla) u -\nu \Delta u +\nabla p =0 \,,\quad \nabla \cdot u=0\,.\label{NS0}
\end{equation}

In a typical engineering or physical  situation  characteristic length, velocity and time scales, $L$ , $U$ and $T$, respectively, are introduced.  After rescaling the relevant quantities in (\ref{NS0}), retaining the same notation for the rescaled dimensionless velocity,  spatial and  temporal  variables, one obtains  the undimensionlized   Navier-Stokes equations

\begin{equation}
\del_t u + (u\cdot \nabla) u - \frac{1}{Re} \Delta u +\nabla p =0 \,,\quad \nabla \cdot u=0\,,\label{NS1}
\end{equation}
 where
\begin{equation*}\Re=\frac{UL}{\nu} \label{Re}\,
\end{equation*}
is the dimensionless Reynolds number, which measures the ratio between a typical  strength of the nonlinear advection term and a typical linear viscous dissipation in the Navier-Stokes equations. Since equations (\ref{NS0}) and (\ref{NS1})  look the same (symbolically),  we will abuse, hereafter, the notation and let $\nu$ represent the  inverse of the Reynolds number.

In most  practical applications the Reynolds number is very large (for instance, the Reynolds number is  about $10^6$  for the air around a moving bicycle,  about  $10^7$  for a moving car, and is of the order of  $10^{12}$ in   meteorological applications). This naturally  suggests to investigate and  compare the behavior of the  Navier-Stokes  equations (\ref{NS1}), for  large values of the Reynolds numbers,   to that of  Euler equations, which are obtained formally by substituting  $\Re=\infty$  in (\ref{NS1}):
\begin{equation}
\del_t u+ (u\cdot \nabla)  u +\nabla p=0\,,\,\, \nabla\cdot u =0\,.\label{Eul}
\end{equation}

Thus, the third thread in our approach would be to compare the behavior of the solutions of (\ref{NS1}) to those of (\ref{Eul}), as the Reynolds number $Re \to \infty$. However,  as it will be developed below, the relationship between these two systems  is far from being simple.

Since the above equations are derived for physical systems they are expected to  respect  certain global conservation laws which indeed can be deduced  by simple computations, {\it but only, assuming that the solutions are smooth enough to legitimize  these computations.}  In  particular, for sufficiently  smooth solutions of the incompressible Navier-Stokes equations, say subject to periodic boundary conditions in a periodic box $\Omega$, by using   the relation:
\begin{equation*}
\int_{\Omega} ((u \cdot \nabla) u)  \cdot u\,dx=\int_{\Omega} \sum_{ij} \del_{x_j} (u_j \frac {u_i^2}2 )\,dx=0 \,,
\end{equation*}
one obtains the formal energy balance:
\begin{equation}
\frac{d}{dt} \int_{\Omega} \frac{|u(x,t)|^2} 2dx +\nu \int_{\Omega} |\nabla u(x,t)|^2dx=0\,.\label{energy0}
\end{equation}

In the incompressible Euler and Navier-Stokes equations the incompressibility condition, $\nabla\cdot u=0$, appears as a constraint on the dynamics of the fluid particles that are governed by the equation (Newton's second law):
\begin{equation}
\frac{d^2 x(t)}{dt^2} =\frac{d}{dt} u(x(t),t) = \nu \Delta u(x(t),t) -\nabla p(x(t),t)\,,
\end{equation}
with the pressure, $p$, being the {\it Lagrange multiplier} corresponding to this constraint. As such it can be eliminated by applying the  Leray projection,   $P_\sigma$,  to the equations (see, e.g., \cite{CF88} for details concerning the Leray projection and other mathematical aspects of Navier-Stokes equations). For a domain without physical  boundaries  (i.e.,  the whole space or periodic box)  $P_\sigma$ is the orthogonal projection from $L^2$ onto the subspace of divergence-free vector fields. That is,  for any given $w \in L^2$,  $P_\sigma w$ is defined by solving the linear system:
\begin{equation}
P_\sigma w=w+\nabla p\,\quad \hbox{and} \quad  -\Delta p= \nabla \cdot w\,,\label{LP1}
\end{equation}
subject to period boundary conditions.
In the presence of physical boundaries the Leray  projection $P_\sigma$  is  the orthogonal projection from $L^2$ onto the subspace of divergence-free vector fields that are tangent to the boundary, hence $P_\sigma$  is defined by solving system (\ref{LP1}), subject to  the Neumann boundary condition
\begin{equation}
\frac{\del p}{\del{\vec n } }=0\,.
\end{equation}

Eventually, due to the fact that our mathematical understanding  of these  problems and of their relations to turbulent flows is far from being complete, the use of explicit examples is important. Such an illuminating  example  is the shear flow solution of Euler equations, which is defined in
$\Omega=\R^3$ or in a periodic box $\Omega=(\R/\Z)^3$ by the formula:
$$ u(x_1,x_2,x_3,t)= (u_1(x_2), 0, u_3(x_1-tu_1(x_2),x_2)\,.$$
Indeed, $u$  is  a solenoidal vector field that solves    Euler equations with the pressure $p=0$.

 Other illustrative examples are the stationary solutions of   Euler equations. Observe that, with $q= p+\frac{|u|^2}2$, Euler equations can also be written as
\begin{equation*}
\del_t u + (\nabla \wedge u)\wedge u + \nabla q=0\,,\quad \nabla \cdot u=0\,.
\end{equation*}
Hence, stationary solutions are characterized by the relation $ (\nabla \wedge u)\wedge u+ \nabla q=0\,.$ In particular, any vector field that satisfies  $\nabla \wedge u= \lambda(x) u$ is a stationary solution with $q=0$ (such solutions  are called Beltrami flows). In particular, any vector field that satisfies   $\nabla \wedge u =0$  forms a stationary solution. In simply connected domains this relation implies that $u$ is a potential flow, i.e., $u=\nabla \phi$, which is  divergence-free, hence $\phi$ is a harmonic potential $-\Delta \phi=0\,.$

Finally, and for convenience we will denote by $(\cdot, \cdot)$ the $L^2$ inner product in the relevant domain.

\section{ The incompressible Euler equations}\label{euler0}

In fact, as indicated by its name, this system of equations  (one of the  first systems of partial differential equations  for fluid mechanics to be introduced) was derived by Euler in 1757, assuming the  incompressibility of the fluid. It has a very rich  mathematical history  with a very interesting  recent important progress.
The Euler equations in a domain $\Omega$ (with or without physical boundaries) are given by

\begin{equation}
\del_t u + (u\cdot\nabla) u +\nabla p=0\, ,\,\, \nabla\cdot u =0\,, \hbox{ in } \Omega \,. \label{Euler1}
\end{equation}
If $\Omega$ has physical boundaries then  the  system (\ref{Euler1}) is supplemented with the no-normal flow boundary condition:
\begin{equation}\label{Euler1-BC}
u\cdot \vec{n}=0 \quad \hbox{on} \quad \partial \Omega\,.
\end{equation}
Observe that the system (\ref{Euler1})
corresponds to the limit case when the viscosity $\nu=0$ in (\ref{NS0}), or the Reynolds number is $\infty$ in   (\ref{NS1}).

In the absence of viscosity (i.e., $\nu=0$)  any {\it regular enough solution } to the Euler system (\ref{Euler1}) conserves the kinetic  energy (compare with (\ref{energy0}) for $\nu=0$) .

It turns out that the vorticity $\omega=\nabla \wedge u$  is ``the basic quantity"  both from the  physical and the mathematical analysis points of view. The system (\ref{Euler1})-(\ref{Euler1-BC}) can be rewritten in an equivalent form in  terms of the vorticity to give
\begin{eqnarray}
&&\del_t \omega + u\cdot \nabla \omega =\omega \cdot \nabla u\,,\label{vort}\\
&&\nabla \cdot u=0\,, \nabla \wedge u=\omega\,, u\cdot \vec n =0\,\hbox{ on } \del\Omega\,. \label{ell2}
\end{eqnarray}
 The elliptic system given in equation (\ref{ell2}) fully determines $u$ in term of $\omega $ (see, e.g., \cite{Ferrari} and references therein), which makes the full  system (\ref{vort})-(\ref{ell2}) ``closed".
 More precisely, the operator $K:\omega \mapsto u$, defined by the unique solution to the elliptic system  (\ref{ell2}), is a linear continuous map from $C^\alpha(\Omega)$ with values in $C^{1,\alpha}(\Omega)$ (with $\alpha>0$), and from the Sobolev $H^s(\Omega)$ with values in $H^{s+1}(\Omega)\,.$

 Furthermore, for $2d$ flows the vorticity is perpendicular to the plane of motion and therefore  equation (\ref{vort}) is reduced (this can also be checked directly) to the advection equation
 \begin{equation}
 \del_t \omega +u\cdot \nabla \omega =0\,.
 \end{equation}
 The quadratic nonlinearity in  (\ref{vort}) has the following consequences, {\it which we describe below, providing only the
 essential arguments, but omitting  the details of the proofs.}

 \subsection{General analytical results for the  $3d$ Euler equations}

 The short time existence and uniqueness of  classical  solutions to the $3d$ incompressible Euler equations, with smooth enough initial data,  was established  (to the best of our knowledge)   by  Lichtenstein in 1925 \cite{LI}  (see also \cite{Majda-Bertozzi}, \cite{MP} and references therein for other proofs).

 \begin{theorem}\label{Lichtenstein}
 Let $u_0 \in C^{1,\alpha}(\Omega)$, with $\nabla\cdot u_0=0$ and $u_0\cdot \vec{n}|_{\del \Omega}=0$. Then there exists a time $T_*>0$, that depends on $u_0$, and  functions $u\in C^0([-T_*,T_*], C^{1,\alpha}(\Omega))\bigcap  \\ C^1([-T_*,T_*], C^{0,\alpha}(\Omega))$, and $p\in C^0([-T_*,T_*], C^{1,\alpha}(\Omega))$, such that the pair $u,p$ solves  Euler equations (\ref{Euler1})-(\ref{Euler1-BC}), with the initial data $u_0$.
\end{theorem}
In the proof of the above theorem, the ``limit" on the lifespan of the  solution is a result of  the following    nonlinear Gr\"onwall estimate, that controls a certain  norm, $y(t)$,  of the solution:
\begin{equation}\label{short-time}
y'\le  C y^{\frac 32} \Rightarrow y(t)\le \frac{y(0)}{(1-2t C y^{\frac12}(0))^2}\,.
\end{equation}
Consequently, for a finite time, which depends  on the size of the initial data $y(0)$, the quantity $y(t)$ is finite. To validate (\ref{short-time}) and the above argument  the  initial data have to be chosen in the appropriate space of regular enough functions. In particular, if we consider the solution in the Sobolev space $H^s$, with $s>\frac  5 2$, then by taking the scalar product, in $H^s$,  of the solution with the Euler  equations (\ref{Euler1}), and using the relevant  Sobolev estimates one obtains:
\begin{eqnarray}\label{scalarprod}
\frac12\frac{d |\!|u |\!|_{H^s}^2}{dt} & =&-((u \cdot \nabla)u , u)_{H^s} =  -(\nabla\cdot (u\otimes u) , u)_{H^s} \\ \nonumber
& = & ( (u\otimes u) : \nabla u)_{H^s}\le C_s |\!|u |\!|_{H^s}^2|\!| \nabla u|\!|_{L^\infty}\le C|\!|u |\!|_{H^s}^3  \,,
\end{eqnarray}
which leads to (\ref{short-time}) and the local in time existence of a smooth solution. Observe that  since $\nabla \cdot u=0$  we have used above the fact that $(u \cdot \nabla)u = \nabla\cdot (u\otimes u)$, where $(u\otimes u)$ is the tensor (matrix) with entries $a_{ij}=u_iu_j$, for $i,j=1,\cdots,d$.

This notion of threshold of regularity for the initial data is well illustrated with the shear flow
$$
u(x,t)=(u_1(x_2),0,u_{3}(x_1-tu_1(x_2),x_2)\,,
$$
which is a well defined solution of Euler equations with pressure $p=0$. Computing the derivative of $u_3$ with respect to $x_2$  yields:
$$
t \del_{x_2}u_1(x_2)\del_{x_1}u_3(x_1-tu_1(x_2),x_2)\,.
$$
 By investigating this term one observes that    initial data which is  in $C^{1,\alpha}$ generate  solutions that will remain in  $C^{1,\alpha}$, for $t\ge 0$,  and this is  in full agreement with  Theorem \ref{Lichtenstein} of Lichtenstein. However,  less regular initial data, say in $C^{0,\alpha}$,  may, in some cases, generate  solutions that will not remain in this space, for $t >0$,  but only in $C^{0,\alpha^2}$  (cf.~\cite{Bardos-Titi}). This counter-example shows that the Euler equations are ill-posed in $C^{0,\alpha}$, i.e. that Theorem  \ref{Lichtenstein}  of Lichtenstein cannot be extended to $C^{0,\alpha}$ initial data.

As it is often the situation  in many nonlinear evolution   problems,  local regularity means:  local existence, uniqueness and local stability (i.e., continuous dependence on initial data). Moreover, one may exhibit a threshold for this existence, uniqueness and propagation of  regularity of the initial data (including analyticity Bardos and Benachour \cite{BB}, see also Kukavica and Vicol \cite{KV} and references therein).
One of the simplest criteria for the regularity, up to the time $T$, is that the stress tensor $\nabla u$ remains bounded in $L^1(0,T; L^\infty)$.  To see that we observe  that this bound implies the integrability of the Lipschitz constant of the vector field $u$, and hence the solutions of the ordinary differential equation   $\frac{dx}{dt}= u(x,t)$ will never intersect and remains smooth within the interval $[0,T]$.  Consequently, and by virtue of (\ref{vort}), the vorticity is transported, along the particles pathes,  and it is stretched by  $\nabla u$, and the latter is controlled by its  $L^1(0,T; L^\infty)$ norm.   However, to underline that it is the dynamics of vorticity which may induce singularities,  a more sophisticated version of this criterion was established   by Beale, Kato and Majda \cite{BKM},   to  only retain the anti-symmetric part of the stress tensor $\big(\nabla u -(\nabla u)^T\big)$ (equivalent to the vorticity vector $\omega$):

\begin{theorem}{\bf Beale-Kato-Majda} \cite{BKM} \label{B-K-M}
Let $u(t)\in C([0,T), H^s) \bigcap C^1([0,T), H^{s-1})$, for $s >\frac{5}{2}$,  be a regular solution of the $3d$ incompressible Euler equations in the interval  $[0,T)$.
Assume that
\begin{equation}
\int_0^T|\!|\nabla \wedge u(.,t)|\!|_{L^\infty} dt <\infty\,,\label{bmk}
\end{equation}
then $u(t)$  can be uniquely extended to the interval  $[0,T+\delta]$, for some  $\delta>0$, as a regular solution of the Euler equations.
\end{theorem}

While a companion result of Constantin, Fefferman and Majda \cite{CFM}  to Theorem \ref{B-K-M}  states that it is mostly the variations in direction of the vorticity that may produce singularities.
\begin{theorem} \cite{CFM}
\label{ Constantin Fefferman Majda}
Let $u$  be a smooth solution of the $3d$ Euler equations, which is defined  in $Q= \Omega\times[0,T)$. Introduce the following quantities (which are well defined for every $t\in [0,T)$):
\begin{equation}
k_1(t)=\sup_{x\in \Omega} |u(x,t)|\,,
\end{equation}
which measures the size of the velocity, and
\begin{equation}
k_2(t)= \sup_{x,y\in\Omega\,,\, x\not=y} \frac{|\xi(x,t)-\xi(y,t)|}{|x-y|}\,,
\end{equation}
which measures the Lipschitz regularity of the direction
$$
\xi(x,t)=\frac{\omega(x,t)}{|\omega(x,t)|}
$$
of the vorticity vector.
Then under the hypotheses:
\begin{equation}
\int_0^T(k_1(t)+k_2(t))dt<\infty\,\, \hbox{and}\,\,\int_0^T k_1(t)k_2(t)dt<\infty
\end{equation}
the solution of the $3d$ Euler equations exists, and is as smooth as the initial data, on the interval $[0,T+\delta]$, for some $\delta>0\,.$
\end{theorem}

It is worth mentioning that a similar criterion to Theorem \ref{B-K-M}, involving the symmetric part of the stress tensor $S(u)=\frac12(\nabla u +(\nabla u)^T)$, was introduced later  by Ponce  \cite{Ponce}. See also \cite{Deng-Hou-Yu} for other related blow-up criteria for the $3d$ Euler equations. Moreover, see  \cite{Gibbon} and references therein for a recent survey of the blow-up problem.

\subsubsection{About the two-dimensional Euler equations}

In the $2d$ case, the vorticity $\omega =\nabla \wedge u $ obeys the evolution equation
\begin{equation}
\del_t \omega  +u\cdot\nabla\omega =0\,.\label{transvort}
\end{equation}
This transport evolution equation enforces  the conservation of any  $L^p$ norm ($1\le p\le \infty$) of the vorticity vector. Taking advantage of this observation in 1963 Youdovitch proved in his remarkable paper   \cite{YU} the global existence and uniqueness  for all solutions with initial vorticity in $L^\infty\,.$ If the vorticity is only  in $L^p$ for $1<p < \infty$ one can prove the global existence (without uniqueness) of weak solutions. The same result holds also for $p=1$, and also for the case when the  vorticity is a finite Radon measure with definite sign (cf. Delort \cite{DE}) or ``simple" changes of sign (cf. \cite{LNX} and \cite{Majda93}; see also \cite{Bardos-Linshiz-Titi}). It is worth mentioning that the proof is more delicate in these limit cases when $p=1$ or when $\omega$ is a bounded Radon measure.

\subsection{Weak solutions of Euler equations} \label{wsol}
\begin{defi}
A  weak solution, on the time interval $[0,T]$, of the Cauchy problem for the  Euler equations (\ref{Euler1}), with divergence-free  initial data $u(x,0)\in (L^2(\R^d))^d$,  is  a divergence-free vector field
$$u(x,t)\in L^2(0,T; R^d)^d\bigcap C([0,T]; (L^2(R^d))_w^d)  \footnote{The symbol$\quad_w$ refers to the weak continuity. In  the present setting this means that for any $\phi\in L^2(R^d)$ the function$t\mapsto (u(t)\, , \, phi(t)) $ is continuous. Moreover such property follows from the relation (\ref{weak0})} $$
and a distribution $p$ which,  for any  smooth test function $\phi(x,t)$ vanishing  for $t$ near $T$, satisfy the Euler equations (\ref{Euler1})  (or  (\ref{Eul}))  in the distribution sense, that is:
\begin{equation}
\begin{aligned}
\int_0^T \int_{\R^d} \left [ \right. u(x,t) \cdot \del_t \phi(x,t)   +& (u(x,t)\otimes u(x,t) ) :\nabla_x \phi(x,t)\left.\right ]\,dx\,dt +\langle p , \nabla_x\phi\rangle  \\
&= -\int_{\R^d} u(x,0) \cdot \phi(x,0)\,dx\,,\label{weak0}
\end{aligned}
\end{equation}
where $\langle\cdot, \cdot\rangle$ denotes the distributional action.
\end{defi}
It is important to observe that the hypothesis $u(x,t)\in L^2(0,T; R^d)^d$  is sufficient to insure that the matrix  $u\otimes u$ is well defined as a distribution (i.e., no need for any additional  hypotheses on the partial  derivatives of $u$).

Moreover, equality   (\ref{weak0}) is equivalent to the following:
\begin{equation}
\begin{aligned}
\int_0^T \int_{\R^d} \left [ \right. u(x,t) \cdot \del_t \phi(x,t)  +&(P_\sigma  (u(x,t)\otimes u(x,t) ):\nabla_x \phi(x,t)) \left.\right ]\,dx\,dt   \\
&= -\int_{\R^d} u(x,0) \cdot \phi(x,0)dx\,,\label{weak1}
\end{aligned}
\end{equation}
with $P_\sigma$ denoting the Leray projection given in (\ref {LP1}). Eventually, this  is also equivalent to assuming that the relation
\begin{equation}
\begin{aligned}
\int_0^T \int_{\R^d} \left [ \right. u(x,t) \cdot \del_t \phi(x,t)  )+& (u(x,t)\otimes u(x,t)) :\nabla_x \phi(x,t) \left. \right]\,dx \, dt \\
&= -\int_{\R^d} u(x,0) \cdot \phi(x,0)\,dx \label{weak2}
\end{aligned}
\end{equation}
holds  for every divergence-free test functions $\phi$.

Moreover,  there are   other variants  of the   notion of weak solutions. One of the most efficient, introduced by DiPerna and Lions \cite{PLLions}, is the notion of {\it dissipative solution}. It is defined  in terms of the stability   with respect to any  smooth divergence-free vector field $w$  (which may not be a solution of the Euler equations). Such a smooth vector field $w$ evidently satisfies the relation:
\begin{equation}\label{Res}
\del_t w+ P_\sigma((w\cdot \nabla) w ) = E \,,
\end{equation}
where $E(w)$ is the residual, satisfying $P_\sigma E=E$. Notice that $E \equiv 0$ if and only if   $w$ is a smooth solution of the Euler equations (\ref{Eul}).
Next, observe that if $u$ is  a smooth divergence-free vector field then by integration by parts, with respect to the spatial variable, one obtains:
\begin{equation*}
\int_\Omega [((u\cdot\nabla) u- (w\cdot\nabla) w) \cdot (u-w)]\,
dx= \int_\Omega (u-w) S(w) (u-w)\, dx \,,
\end{equation*}
with
\[
S(w)=\frac12(\nabla w+(\nabla w)^T)\,.
\]
If in addition we assume  that $u$ is a smooth solution of the Euler equations then from (\ref{Res}), (\ref{Eul}) and the relevant boundary conditions (e.g., periodic boundary conditions or  (\ref{Euler1-BC})) we conclude :
\begin{equation*}
\begin{aligned}
\frac{1}{2}\frac{d}{dt}\int_\Omega|u-w|^2\, dx=  \int_\Omega (u-w) S(w) (u-w)\, dx +\\
 \int_\Omega E(w(x,t)) \cdot (u(x,t)-w(x,t))\,dx \,,
\end{aligned}
\end{equation*}
which, upon integration with respect to time, gives
\begin{equation*}
\begin{aligned}
&\int_\Omega |u(x,t)-w(x,t)|^2\le \\
&2 \int_0^t\int_\Omega |E(w(x,s)) \cdot (u(x,s)-w(x,s))|\,dx\,ds
\\
&+2\int_0^t \int_\Omega |( u(x,s)-w(x,s)) S(w) (u(x,s)-w(x,s))|\,dx \, ds\\
&+ \int_\Omega |u(x,0)-w(x,0)|^2\,dx\,;
\end{aligned}
\end{equation*}
and by virtue of  Gr\"onwall inequality yields
\begin{equation}
\begin{aligned}
&\int_\Omega |u(x,t)-w(x,t)|^2 \,dx\le  \Big (\int_\Omega |u(x,0)-w(x,0)|^2\,dx \Big ) e^{2\int_0^t|\!|S(w)|\!|_\infty(s) ds} \\
&+2 \int_0^t \Big (\int_\Omega |E(w(x,s)) \cdot (u(x,s)-w(x,s))|\,dx \Big) e^{2 \int_s^t|\!|S(w)|\!|_\infty(\tau) d\tau}\,ds.\label{dissol}
\end{aligned}
\end{equation}
Inspired by the above we introduce the following definition of dissipative solution for  Euler equations (\ref{Eul}):
 \begin{defi}
 A  dissipative solution of  Euler equations (\ref{Eul}) is a divergence-free vector field $u$ which together with  any smooth divergence-free vector field $w$ (which  need not be a solution of  Euler equations) they satisfy the inequality:
\begin{equation}
\begin{aligned}
&\int_\Omega |u(x,t)-w(x,t)|^2\, dx \le  \Big(\int_\Omega |u(x,0)-w(x,0)|^2\, dx\Big)\, e^{2\int_0^t|\!|S(w)|\!|_\infty(s) ds} \\
&+2 \int_0^t|(\del_t w(s)+P_\sigma ((w(s)\cdot \nabla) w(s)) , u(s)-w(s))|\,e^{2 \int_s^t|\!|S(w)|\!|_\infty(\tau) d\tau}\,ds\,.\label{dissol2}
\end{aligned}
\end{equation}
 \end{defi}

\begin{remark}\label{disssolr}
With the above definition one should keep in mind the following facts.

\noindent 1. Inequality  (\ref{dissol2}) makes sense with no other regularity conditions on $u$ other than the fact that it belongs to $(L^2([0,T;(L^2(\Omega))^d))^d$ ($d=2$ or $d=3$).

\noindent 2.  If $w$ is a genuine smooth solution to Euler equations (\ref{Eul}) then one has $E(w)\equiv 0$ and inequality (\ref{dissol2}) reduces to
\begin{equation}
\int_\Omega |u(x,t)-w(x,t)|^2\, dx\le  \Big(\int_\Omega |u(x,0)-w(x,0)|^2 \, dx\Big) \,e^{2\int_0^t|\!|S(w)|\!|_\infty(s) ds} \,,
\end{equation}
which implies both uniqueness of dissipative solutions, within the class of smooth solutions, and  stability  of dissipative solutions with respect to smooth solutions.
\end{remark}

Below it will be observed that weak solutions of Euler equations (\ref{Eul})  may not conserve energy. Therefore, the connection between weak and dissipative solutions of (\ref{Eul}) requires some additional information on the control of  energy.
Following \cite{CLLS1} we introduce the following:
\begin{defi}\label{admissible-def}
A divergence-free vector field $u\in L^2(0,T; R^d)^d$ is called an admissible weak solution if it is a weak solution and if it satisfies the additional hypothesis:
\begin{equation}
\int_\Omega |u(x,t)|^2\, dx \le \int_\Omega |u(x,0)|^2 \, dx\,,\qquad \hbox{for all} \quad  t \in [0,T]\, .\label{admiss0}
\end{equation}
\end{defi}
As a result of the above definitions we have:
\begin{proposition}
 Let $u$ be an  admissible weak solution of  Euler equations (\ref{Eul}), in $\Omega=\R^d$ or in a periodic domain $\Omega=(\R/\Z)^d$,  then $u$  is a dissipative solution.
\end{proposition}
{\bf Proof}
 Below, we only present  the main idea of the  proof.  Let $w$ be any smooth divergence-free vector field. The point is to compare
$$
\frac12\int_\Omega |u(x,t)-w(x,t)|^2dx =\frac12\int_\Omega |u(x,t)|^2dx + \frac12\int_\Omega |w(x,t)|^2dx -\int_\Omega (u(x,t) \cdot w(x,t))dx
$$
and
$$
\frac12\int_\Omega |u(x,0)-w(x,0)|^2dx =\frac12\int_\Omega |u(x,0)|^2dx   +  \frac12\int_\Omega |w(x,0)|^2dx -
\int_\Omega(u(x,0) \cdot w(x,0))dx
$$
With hypothesis  (\ref{admiss0} ) it is then enough to show that:
$$
\begin{aligned}
&\frac12\int|w(x,t)|^2-\int(u(x,t) \cdot w(x,t)) \, dx=\frac12\int|w(x,0)|^2-\int(u(x,0) \cdot w(x,0)) \, dx\\
&+\int_0^t\frac{d}{ds}\left (  \frac12\int|w(x,s)|^2-\int(u(x,s) \cdot w(x,s))\, dx\right )ds\\
&\le \frac12\int|w(x,0)|^2-\int(u(x,0) \cdot w(x,0)) \, dx\\
& {+\int_0^t\int_\Omega (u-w) S(w) (u-w)\, dxds +}\\
& { +\int_0^t\int_\Omega E(w(x,s)) \cdot (u(x,s)-w(x,s))\,dxds \,,}
\end{aligned}
$$
or that
\begin{equation*}
\begin{aligned}
&\frac{d}{ds} \left (\frac12\int|w(x,s)|^2-\int(u(x,s)\cdot w(x,s))dx\right )= {\int_\Omega (u-w) S(w) (u-w)\, dx +}\\
& { \int_\Omega E(w(x,s)) \cdot (u(x,s)-w(x,s))\,dx \,.}
\end{aligned}
\end{equation*}
Here, the most   nontrivial step that needs further justification is the one involving the term
$\int_\Omega(\del_s u(x,s) \cdot w(x,s))\, dx$.  But   since $w$ is smooth vector field  the notion of weak solution can be invoked.

However, for Euler equations the introduction   of weak solutions corresponds to an important loss of information illustrated by the following theorems of De Lellis, Sz\'ekelyhidi and Wiedemann.
\begin{theorem}\label{DSW}
\noindent 1.  Let $\Omega\subset \R^d$, for $d=2,3$,  be any  bounded open set, and let $T >0$.  Then for any  $e(x,t)>0$, continuous function defined on $Q=\overline {\Omega \times[0,T]} \subset \R^d\times\R_t$,  there exists an infinite set of weak solutions of   Euler equations in $\R^d$ which satisfy:
\begin{equation*}
 \frac{|u(x,t)|^2}2  =e(x,t), \,\hbox{for all} \, (x,t) \in Q, \quad \hbox{and } \,|u(x,t)| =0,\,   \hbox{for all} \,  (x,t)\notin Q .
\end{equation*}

\noindent 2.  There is an infinite set of initial data $u_0\in L^2(\R^d)$  for which there exist admissible solutions (hence they are dissipative).

\noindent 3. For any initial data $u_0\in L^2(\R^d)$ there exists infinitely many  weak solutions  $u(x,t) \in C(\R_t; L^2(\R^d)_w)$, such that $u(x,0)=u_0$.

\end{theorem}

The first two statements in Theorem \ref{DSW}  are due to De Lellis and Sz\'ekelyhidi, and the last one is due to Wiedemann following their approach. A complete proof of the above theorem is out of the scope of this contribution, but in the following subsection we give some comments which might be appropriate.

{\bf Comments on Theorem \ref{DSW} }  The first two statements of  Theorem \ref{DSW}   have, in fact, a long history over the second half of the previous century. In 1993 Scheffer  \cite{SC} was  first person to observe the existence of weak solutions of Euler equations with compact support in time. Later, in 1997,   Shnirelman \cite{Sh} refined  Scheffer's proof. The basic method was to construct approximate highly oscillating (in space and time) divergence-free vector fields $u_k$ (with compact support in space and time),  that would converge strongly in $L^2(\R^d\times \R_t)$  (in order to be able to  pass to the limit in the nonlinear term $u_k\otimes u_k$), while by equation (\ref{Res})  the sequence of residuals:
$$
E_k:= \del_t u_k + P_\sigma((\nabla \cdot ( u_k\otimes u_k))\,,
$$
would converge weakly to zero.

Of course, such weak solutions, for which De Lellis and Sz\'ekelyhidi  coined the name ``wild", are  not   physical at all. They describe  a fluid flow which starting   from rest and without any external force applied to it, begins to  move and then goes back to rest. If this situation could be applied to ``real world fluid flows " that would solve the energy crisis!

To establish their proof,   De Lellis and Sz\'ekelyhidi   realized a striking analogy between   ``wild solutions" of Euler equations  and the classical problem of isometric imbedding.

The simplest example of such an imbedding problem goes as follows:

 Let $S^1 =\{ x\in \R^2:\,|x|=1\}\subset \R^2$  be  the unit circle, and,   for any $r>0$,   $B^2(r)=\{ x\in \R^2: \,|x|\le r\}\subset \R^2$ be the closed disc of radius $r$. The question is: what is  the best regularity (in term of the H\"older exponent $\beta$)  that is independent of $r$, and  for which there exists a map $\phi_r\in C^\beta(S^1;B^2(r))$  (dubbed isometric imbedding)  which enjoys  the following property:
  $\phi_r(S^1)$ is, for the Riemannian metric induced by $\R^2$, a manifold isometric to $S^1$?

%
Now if one requires $\phi_r \in  C^2$, the curvature of $\phi_r(S^1)$ should be bounded from above and below, independent of $r$. While to place the curve  $\phi_r(S^1)$ in a ball of arbitrarily small radius, one should expect its curvature to go to infinity, as $r\to 0$. This is an obstruction which was realized by Gauss under the name of ``egregious  theorem". Hence, to realize the isometric imbedding, one has to relax the regularity of the mapping $\phi_r$ from being $C^2$.   What was proven (cf. \cite{CLLS2} and the reference therein) is the existence of  a small enough H\"older exponent $\alpha\in (0,1)$, which is $r$ independent,  such that for any $r>0$ there exists an isometric imbedding $\phi_r \in C^{1,\alpha}$ of $S^1$ into $\phi_r(S^1) \subset B^2(r)$.
The main idea of the proof is to show first the existence of a mapping $\psi_r$ that reduces the distance, on the Riemannian manifold $\psi_r(S^1)$, and then  accumulate (add)  small oscillations (wiggles of the curve) to reach the equality of the lengths; and thus prove that the isometric imbedding can be realized with the limit function of this process, $\phi_r \in C^{1,\alpha}$.

In the same way one can rewrite  Euler equations as a linear system with a nonlinear constraint,  namely,  it is equivalent to the system:
\begin{equation}
\del_t u + \nabla_x \cdot M + \nabla q=0 \,,\quad u(x,0)=u_0(x)\,,
\nabla \cdot u=0\,,\label{tartar1}
\end{equation}
with the constraint
\begin{equation}
M=u\otimes u-\frac{|u|^2}d I\,,\label{tartar3}
\end{equation}
where $M$ is a traceless self-adjoint matrix, and $I$ is the identity matrix. To solve the above system one considers, instead, the following more relaxed problem. Let $Q$ and $e$ be as in Theorem \ref{DSW}.
Let $\Lambda$ be is the closed hull  (in some convenient topology) of the set of all  triples $(u,q,M)$, where $u$ is a vector field,   $q$ is a scalar, and $M$ is a self-adjoint traceless matrix, which satisfy  system (\ref{tartar1}),  and the following relaxed constraint:
\begin{equation*}
 (u\otimes  u)-M \le \frac{2 e(x,t)}d I,\, \hbox{for all}\, (x,t) \in Q,\, \hbox{ and}\,\,  |u(x,t)| =0\,, \hbox{for all} \, (x,t)\notin Q. 
\end{equation*}
The elements of  $\Lambda$ are called sub-solutions of (\ref{tartar1})-(\ref{tartar3}). Then the proof of Theorem \ref{DSW} goes as follows. First, one should  show that the set of sub-solutions, $\Lambda$,  is not empty. Then, to conclude the proof,  one constructs a sequence $(u_k,q_k,M_k)\in \Lambda$, which is constructed by adding, more and more, spatial and temporal   oscillations of small amplitudes  and high  frequencies, and   which converges strongly, in $L^2$, such that:
\begin{equation}
 \lim_{k\rightarrow \infty}\left [ (u_k\otimes  u_k)-M_k - \frac{|u_k|^2}d I \right ]=0\,.
\end{equation}

\begin{remark}
1. Uniqueness among dissipative solutions,  whenever one of them is regular, implies that the wild  solutions stated in  point 3 of Theorem \ref{DSW}  (for every initial data $u_0\in L^2$) cannot be admissible whenever the initial data is regular. Indeed, as Wiedemann  observes in \cite{Wie} his  construction of sub-solutions  requires a jump in energy   at time $0^+$.

\noindent 2. Similarly the initial data giving rise in  point 2 of Theorem \ref{DSW} to admissible solutions cannot be regular; otherwise,  the solutions being admissible  would coincide with the unique dissipative solution. On the other hand, since the proof involves the Baire  category theorem,  the construction of these initial data is,  in general, far from being explicit. However,  for the case of  Kelvin-Helmholtz type initial data:
\begin{equation}
u_0(x_1,x_2) =(-\sgn(x_2),0)\,,
\end{equation}
 in a periodic domain, explicit sub-solutions, and then infinitely many solutions, have been provided in \cite{vortexpaper}. In   $3d$, and  for initial data of the form  $(-\sgn(x_2),0,w(x_1,x_3)$, the   construction presented in \cite{vortexpaper} gives  rise to  an infinite set of admissible solutions, cf. \cite{BTW}.

\noindent 3. Eventually, with a better adapted choice of the high-frequency corrections to the sub-solutions, one can improve the spatial and temporal regularity of the wild solutions obtained at the limit. Along these lines, using Beltrami flows as building blocks for the high-frequency corrections in \cite{CLLS3}, one proves, for every fixed $\alpha$ in $(0,1/10)$, the existence of admissible wild solutions which belong to $C^{0,\alpha}$,  and which exhibit strict   energy decay.
\end{remark}

\subsection{Energy decay and turbulence}

For   then Euler equations the only natural estimate is the energy estimate which, as mentioned above, gives for smooth solution the relation:
\begin{equation}
\frac12\int|u(x,t)|^2dx =\frac12\int|u(x,0)|^2dx\,.
\end{equation}
From this relation one deduces that any $\epsilon $ dependent family $u_\epsilon$ of smooth solutions
 is  uniformly bounded according to the formula:
 \begin{equation}
 \forall t ,\epsilon\quad |u_\epsilon(t)|_{L^2} \le \sup_{\epsilon} |u_\epsilon(0)|_{L^2}<\infty\,.
 \end{equation}
 Variants  of  this energy estimate (with inequality instead of  equality) can be obtained and are valid in most of the physical situations (including Navier-Stokes  with weakly converging initial data  and viscosity going to $0$ together with boundary effects) .

 This implies the existence of subsequences, still denoted by $u_\epsilon$, which are  weakly converging to a limit $(\lim_{\epsilon \rightarrow 0} u_\epsilon ) (x,t)\in L^{\infty} (\R_t;L^2)\,.$ However,  at this level of generality it is not possible to conclude that this limit remains a solution of Euler equations. The problem is that the weak limit of the tensor product may not be the product of the weak limits:
 \begin{equation}
 \lim_{\epsilon \rightarrow } (u_\epsilon \otimes u_\epsilon)\not= (\lim_{\epsilon \rightarrow 0} u_\epsilon ) \otimes( \lim_{\epsilon \rightarrow 0} u_\epsilon )
 \end{equation}
and therefore the weak limit $\overline u=\lim_{\epsilon \rightarrow 0} u_\epsilon  $ is solution of the equation:
\begin{equation}
\begin{aligned}
&\nabla \cdot \overline u =0\,, \quad \del_t \overline u + \nabla(\overline u \otimes \overline u) +  RT + \nabla \overline p =0\,,\\\label{rt}
&\hbox{ with } \quad RT=\lim(u_\epsilon -\overline u)\otimes \lim(u_\epsilon -\overline u)\,.
\end{aligned}
\end{equation}
One simple example (which is due to   DiPerna and Majda \cite{DM}) can be constructed with the shear flow:

\begin{equation*}
u_\epsilon(x,t)= ( \sin({\frac{x_2}\epsilon}),0,  w(x_1-t\sin({\frac{x_2}\epsilon}))
\end{equation*}
which is uniformly bounded and which weakly converges to
$$
\overline u=(0,0,\int_0^1w(x_1-t\sin s)ds)\,.
$$  Obviously this  is not a solution of   the Euler equations, by direct inspection and also because the term
$u^1_\epsilon\del_{x_1} w_\epsilon $ does not go to zero!

Weak limits can be interpreted in term of ``coarse graining", as the deterministic counter-part of the Reynolds stress tensor for small fluctuations in the statistical  theory of turbulence. In this setting one could consider, as turbulent flows which  are not solutions of  the Euler equations due to the appearance  of the extra term
$RT$ in    (\ref{rt}).
Now  the formula
\begin{equation}
\hbox{ for almost every} \, t \quad \lim_{\epsilon\rightarrow 0} |\overline u(t) -u_\epsilon(t)|^2 =\lim_{\epsilon\rightarrow 0} |u_\epsilon(t)|^2-|\overline u(t)|^2
\end{equation}
 implies that the conservation of   energy would give $RT=0$ and  therefore  no turbulence. Hence, this implies the  conjecture that turbulence would be related to anomalous energy dissipation.

 Indeed, the issue of the conservation of energy for weak solutions  was already raised  by Onsager \cite{ON} when he  formally wrote:
\begin{equation}
0= \int [\nabla\cdot(u\otimes u)]\cdot udx\simeq  \int (\nabla ^{\frac13} u \cdot \nabla ^{\frac13}u\cdot   \nabla ^{\frac13}u) \simeq \int | \nabla ^{\frac13} u|^3dx \label{naive}\,.
\end{equation}
He conjectured that any weak solution that would be more regular than $C^{0,\frac13}$ would preserve  energy. A rigorous proof of this result, inspired by the naive formula (\ref{naive}), was given by Constantin, E and Titi \cite{CET}.
Of course this does not mean that the conservation of energy implies the regularity of the solution; simple examples of singular weak solutions which preserve energy can be constructed  with the shear flow \cite{Bardos-Titi} or, as noticed  above following (\cite{CLLS2}), as ``wild" solutions.
Eventually to complete this picture it would be interesting to exhibit weak solutions less regular than $C^{0,\frac13}$ for which  energy would decay.

The first  step was the construction of  continuous (in space and time) solutions with decaying energy (cf.  \cite{CLLS4}  for the spatially three-dimensional case,    and  \cite{CCLLS4} for the spatially two-dimensional case.)  One of the basic ingredients in the spatially three-dimensional case  was the introduction of the Beltrami flows in the construction.   Then a refinement of the above method
gave the existence of energy decaying weak solutions  which belong to $C^{0,\alpha}$,  with $\alpha< \frac1{10}$  (cf. \cite{CLLS3} for the spatially three-dimensioanl case, and  \cite{choffrut} for the spatially two-dimensional case). Most recent improvement, due to Isett \cite{Isett}, extend the H\"older regularity of these energy decaying solutions to $\alpha <\frac15$,  in a preprint of 171 pages!

\section{Navier-Stokes equations without and with  boundary effects}

\subsection{Why in some cases Navier-Stokes rather than Euler} In the section \ref{preliminaries} the Euler equations were presented as a formal  limit, as $Re\rightarrow \infty$, of Navier-Stokes equations. However, we have also added the the relation between these two equations   ``is far from being simple."    The Euler equations were  derived first, and it
  took about 45  years to obtain   Navier-Stokes equations as an  ``improved"  model  (cf. Navier \cite{NAV} and Stokes, \cite{STO}), taking  into account the viscosity effects which  are better adapted to the real fluids . Of course, at that time this was not motivated by the arguments of De Lellis and Skz\'ekelyhidi on wild  solutions,   but by the need for a mathematical description of the fluid's viscosity. Moreover, the introduction of   the viscous term $-\nu \Delta u$,  with the corresponding boundary conditions,  is essential to describe the interaction of  a fluid flow with physical obstacles. In general, this physically observed interaction  cannot be described only by means of  Euler equations, subject to the impermeability (no-normal flow) boundary conditions, as it will be illustrated by the following proposition:
  \begin{proposition}
  Consider an obstacle $K\subset \R^3$, with a smooth boundary which is topologically isomorphic to a closed ball. Let   $u$ be the velocity field of a stationary solution of Euler equations in $\Omega= \R^3\backslash K$, that is constant  at infinity, i.e.,  $\lim_{|x|\rightarrow \infty } u(x)=u_\infty $,  and which is satisfying the impermeability (no-normal flow) boundary condition, $u\cdot \vec n =0\,$ on $\del \Omega$. Then  the total force exerted by the fluid flow on the obstacle $K$ is $0$.
  \end{proposition}

  {\bf Proof } (cf.~\cite{CM}  or \cite{MP}  for details)
  Since the fluid flow  is  stationary and constant at infinity it is also irrotational and therefore one has the Bernoulli  theorem:
  \begin{equation}
  \nabla \Bigg(\frac {|u|^2}2 + p\Bigg)=0 \label{bernoulli}\,.
  \end{equation}
  Moreover, since the obstacle $K$ is homeomorphic to a ball in $\R^3$  the domain $\Omega = \R^3 \backslash K$ is simply connected; therefore the velocity field $u$ fluid is both ``potential" and solenoidal .  Hence, it can be described by the system:
  \begin{equation}
  u= -\nabla \phi\,,\quad -\Delta \phi =0\label{lap0}
  \end{equation}
  Thus,  with the asymptotic behavior, at infinity, of the Green function of the Laplacian in $\R^3$ (behaving like $(4\pi)^{-1} |x|^{-1}$) one has:
  \begin{equation}
  u(x)  = u_\infty + O(|x|^{-3})\, \hbox{and } p(x)=p_\infty + O(|x|^{-3})\,,\label{green}
  \end{equation}
  where the above asymptotic behavior of the pressure follows from (\ref{bernoulli}).
  Now, to compute the force exerted by the fluid one introduces  the ball $B_R=\{ x\in \R^3, |x|<R\}$. For $R$ large enough one has $K\subset B_R$. Therefore, thanks to the  impermeability boundary condition, $u\cdot\vec n=0$,   one  can compute, $\vec F$,
  the force exerted on $K$, as follow:
 \begin{equation}
\begin{aligned}
\vec F= &-\int_{\del\Omega}p \vec n d\sigma
= -\int_{\del\Omega}\Bigg((u\cdot \vec n) u + p \vec n\Bigg) d\sigma \\
&+\int_{|x|=R}\Bigg((u\cdot \frac{x}{R} ) u + p \frac{x}{R}\Bigg) d\sigma
-\int_{|x|=R}\Bigg((u\cdot \frac{x}{R} ) u + p \frac{x}{R}\Bigg) d\sigma\\
=& - \int_{\Omega \cap B_R}\Bigg( (u\cdot \nabla) u +\nabla  p\Bigg)dx  -\int_{|x|=R}\Bigg((u\cdot \frac{x}{R} ) u + p \frac{x}{R}\Bigg) d\sigma \\
=& - \int_{|x|=R}\Bigg((u\cdot \frac{x}{R} ) u + p \frac{x}{R}\Bigg) d\sigma\label{dal} \,,
\end{aligned}
\end{equation}
where the last equality is due to (\ref{bernoulli}).
Finally, by virtue of  (\ref{green}) the relation  (\ref{dal}) implies:
\begin{equation}
\vec F= -\int_{\del\Omega}p \vec n d\sigma =-\lim_{R\rightarrow \infty }\int_{|x|=R}\Bigg((u\cdot \frac{x}{R} ) u + p \frac{x}{R}\Bigg) d\sigma=0\,!
\end{equation}
\begin{remark}
The above computation was already known by d'Alembert, for this reason  it is called the d'Alembert paradox.
\noindent
 A similar computation can be performed in two dimensions, but   in this case the domain $\Omega = \R^2 \backslash K $, where $K$ is isomorphic to a disc in $\R^2$ is not simply connected. Therefore, irrotational vector fields need not to be potential. In addition,  the Green function in two dimensions behaves like $\frac1{2\pi} \log \frac1{|x|}$, and  it produces a nonzero contribution (lift) called the Kutta-Joukowski relation.

\end{remark}
 In conclusion, we observe that the introduction of the viscosity, and of  the corresponding  boundary condition associated with  to the Laplacian,  may lead to effects that cannot be described by Euler equations (which corresponds to the singular limit $\nu=0$),  and that these effects might become important for  ``very small" values of $\nu$, where $\nu>0$.

 For example, the $2d$ Navier-Stokes equations in the half-plane, $x_2>0$, subject to Dirichlet boundary condition, is a well-posed problem; and its vorticity, $\omega$, evolves according to the equation
 \begin{equation}\label{Vor-2D}
 \del_t \omega_\nu -\nu \Delta \omega_\nu +u_\nu \cdot \nabla \omega_\nu =0\,.
  \end{equation}
  In the inviscid case, $\nu=0$, the vorticity, $\omega$, evolves according to:
 \begin{equation}
 \del_t \omega  +( u \cdot \nabla) \omega=0\,,
 \end{equation}
  which  does not indicate any interaction of the vorticity inside the domain with the boundary  $\R \times \{x_2=0\}$.  However, for the viscous case, $\nu >0$,   the values of the vorticity at the boundary, $\R \times \{x_2=0\}$,  will depend in  a nontrivial way on the pressure .  To see that, we  compute the normal derivative of the vorticity at $\R \times \{x_2=0\}$, and use the incompressibility condition,  $\nabla\cdot u_\nu=0$, and the Dirichlet boundary condition, $u_\nu =0 $ at $\R \times \{x_2=0\}$,  to obtain:
$$\partial_{x_2} \omega_\nu = \partial_{x_2}(\partial_{x_1} u_\nu^2-\partial_{x_2} u_\nu^1)=
\partial_{x_1}(\partial_{x_2} u_\nu^2)-\partial_{x_2}^2u_\nu^1=
-\partial_{x_1}^2 u_\nu^1-\partial_{x_2}^2u_\nu^1=-\partial_{x_2}^2u_\nu^1 \,.
$$
Next, use the evolution  equation for  the first component of the velocity, $u_\nu^1$, in the Navier-Stokes equations at the boundary, $\R \times \{x_2=0\}$,  to get:
$$
\nu \del_{x_2}^2u_1= \del_{x_1} p\,.
$$
Combining the above yields the so called Lighthill equation:
\begin{equation} \label{Vor-BC}
-\nu \del_{x_2} \omega_\nu=\del_{x_1} p, \, \hbox{at} \,  \R \times \{x_2=0\}\,.
\end{equation}
The above relation, (\ref{Vor-BC}),  implies that the variation of pressure at the boundary generates the vorticity there  or vice versa. Notice that the quantities involved in  (\ref{Vor-BC}) may be very large, for small values of $\nu >0$,  leading (cf. section (\ref{vbdef})) to important effects and interaction with the boundary.

To measure such effects, in particular in comparison with the Reynolds number, it is convenient to consider a slightly more general class of boundary conditions. Therefore, in  an open smooth domain  $\Omega\subset \R^d$, with $d=2$ or $d=3$, one supplements the Navier-Stokes equations
\begin{equation}
 \del_t u_\nu + (u_\nu\cdot \nabla) u_\nu -\nu \Delta u_\nu+ \nabla p_\nu=0 \,,  \label{genns}\hbox{ in } \Omega
 \end{equation}
with  the boundary conditions:
\begin{eqnarray}
&&{u_\nu \cdot \vec n =0\,,\quad  \nu ( \del_{\vec n}u_\nu  + C(x) u_\nu )_\tau+\lambda u_\nu =0\,\hbox{ on } \del \Omega\,,  \label{boundcond}}\\
&&\hbox{ with }\quad \lambda(\nu,x)\ge 0 \,,\,\, C(x)\in C( \R^n; \R^n) \,,\label{coeff}
\end{eqnarray}
 where $\vec n$ is denoting the outward normal to $\del \Omega$, and $\vec{a}_\tau = \vec{a} - (\vec{a}\cdot \vec{n}) \vec{n}$ is the tangential component of the vector $\vec{a}$ at the boundary $\del \Omega$.  Of course,  when the problem is considered in the whole space, or in a periodic box, these boundary conditions are  omitted. In (\ref{boundcond}) there are two statements which represent,  on one hand, the impermeability condition $u_\nu\cdot \vec n=0\,,$ which is in some sense universal (independent of the fluid, and in particular of whether it is viscous or inviscid)
and,  on the other hand, a condition which models the friction between the fluid flow and the boundary, and which may exhibit some nontrivial behavior in terms of the Reynolds number.
With $\lambda (\nu)=\infty$, the second statement corresponds to the Dirichlet boundary condition:  the component of the velocity on the plane tangent  to $\del \Omega$ being also zero and, therefore, it is called the ``no-slip" boundary condition.

In the more general case, $C(x)$ is  a linear bounded operator acting on the space of tangent vectors to $\del\Omega$ with the $L^2(\del\Omega)-$norm.
With a convenient choice of $C(x)$ this statement includes the following conditions
\begin{equation}
(S(u_\nu)\cdot \vec n)_\tau = (\del_{\vec n}u_\nu)_\tau -(\nabla^t \vec n \cdot u_\nu)_\tau\Rightarrow \nu(S(u_\nu)\cdot \vec n)_\tau +\lambda u_\nu=0\,,\label{NF}
\end{equation}
and
\begin{equation}
(\nabla \wedge u_\nu)\wedge \vec n=(\del_{\vec n}u_\nu)_\tau +(\nabla^t \vec n \cdot u_\nu)_\tau\Rightarrow \nu(\nabla \wedge u_\nu)\wedge \vec n +\lambda u_\nu=0\,.
\label{Vort}\end{equation}
Condition (\ref{NF}) is  the standard Fourier law, which can be either derived by  phenomenological arguments, or deduced from a boundary condition for the Boltzmann equations, due to the atomic structure of the wall (cf. \cite{AIO}, \cite{BGP} and \cite{MLRS}).  Condition (\ref{Vort})  emphasizes  the role of the vorticity.

Integrating by parts after multiplying equation (\ref{genns}) by $u_\nu$, and using some standard trace estimates one obtains the:
\begin{proposition}
Any smooth solution of the Navier-Stokes equations, (\ref{genns}), with the boundary conditions (\ref{boundcond}), satisfies the relation (energy balance):
\begin{equation}
\begin{aligned}
\frac 12\frac{d}{dt} \int_\Omega |u_\nu(x,t)|^2dx &+ \nu \int_\Omega |\nabla u_\nu|^2dx \\
&+\int_{\del \Omega} \lambda(x)|u_\nu(x,t)|^2d\sigma
= \nu \int_{\del \Omega} C(u_\nu)_\tau u_\nu d\sigma\,;\label{enform}
\end{aligned}
\end{equation}
and the $\nu-$uniform estimate:
\begin{equation}
\frac{d}{dt} \int_\Omega \frac{|u_\nu(x,t)|^2}2dx + \nu \int_\Omega |\nabla u_\nu|^2dx +\int_{\del \Omega} \lambda(x)|u_\nu(x,t)|^2d\sigma
= o(\nu)\,.\label{enform2}
\end{equation}
\end{proposition}

\subsection{Basic mathematical results for the Navier- Stokes equations}

From   the estimate (\ref{enform}) one deduces, and this was the original contribution of Leray  \cite{Leray},  that for any initial data $u_0\in L^2 $ there exists (globally in time) weak solution to the Cauchy problem. However,  there are several important differences between the $2d$ and $3d$ cases.  In the $2d$ case  one can prove that for any $u_0\in L^2$ the solutions are unique and that the formal energy equality (\ref{enform}) is exact and can be justified rigorously.
In the $3d$ case the situation (which is an object of the present Clay prize) is more complex. With no restriction on the size of the viscosity, $\nu$, when it is compared inversely to the spatial oscillation in the initial data (even very smooth initial data), the only result available is the (global in time) existence of a weak solution which satisfies the energy inequality, i.e. instead of the  equality sign $=$  in (\ref{enform}) one has the inequality sign $\leq$. Such weak solutions that satisfy the energy inequality in (\ref{enform}) are called Leray-Hopf solutions. It is worth mentioning that all the known existence proofs of weak solutions to the Navier-Stokes equations lead at the end to Leray-Hopf solutions. The uniqueness of weak solutions, or Leray-Hopf weak solutions, and the propagation of regularity, for any time, of the initial data are still  open problems.

There  is a long history of question of global regularity and of blow up criteria  of smooth solutions of the $3d$ Navier-Stokes equations, taking a full advantage of    the regularizing effect  of the Laplacian.
As observed in Theorem \ref{B-K-M} the absence of blow up for the solutions of the $3d$ Euler equations is enforced by the Beale-Kato-Majda condition, which concerns the vorticity. For the Navier-Stokes equations it is enough to assume that  the velocity remains bounded to ensure the regularity. \footnote{The appearance of unbounded, or very large velocity, is not in inconsistent  with the physical derivation of the equations. This is because, as it has been mentioned in   the ``Preliminaries" section (section \ref{preliminaries}), the incompressible Navier-Stoke equations can be deduced as a small Mach number approximation. In this case the velocity $u_\epsilon$ is replaced by $\epsilon \tilde u_\epsilon$, and the fact that $u_\epsilon$ is bounded does not imply that $\tilde u_\epsilon$, the solution of the ``incompressible" equations, needs to be also uniformly bounded.}  This condition can be widely improved starting from the pioneer works of Prodi \cite{Prodi} and Serrin \cite{Serrin}, involving the $L^p((0,T), L^q(\Omega))$ norms of the velocity field, with $\frac1 p + \frac 1 q = 1$ (for more details see, for example,  the survey articles \cite{Germain} \cite{LADY03}, \cite{TT00} and references therein). Most recently, there has also been some criteria for  global regularity of the $3d$ Navier--Stokes equations involving  the pressure (cf.~\cite{CT08} and references therein). In addition, some other sufficient regularity conditions were recently
established in terms of only one component of the stress tensor $\nabla u$  (cf. \cite{Cao-Titi} and references therein). Furthermore, in the presence of certain spatial symmetries, such as axi-symmetric \cite{LadyPaper} or helical \cite{MTL}  flows, the Navier-Stokes equations are known to be globally well-posed, but such results are still out of reach for Euler equations (see, e.g., \cite{Ettinger} and references therein for some related results concerning helical Euler equations).  It is worth mentioning that  vorticity stretching term, which is the main characteristic of $3d$ flows verses $2d$ flows, is nontrivial for  axi-symmetric and  helical   flows.

Eventually,   the two following sharp results  (even, if once someone proves that weak solutions of Navier-Stokes equations are always smooth)  have a genuine interest by the detail and the physical interpretation of their proofs.

\bu$\,$ Caffarelli, Kohn and Nirenberg \cite{CKN}:

For any initial data $u_0\in L^2$ there is a weak solution (called   a suitable solution) for which the support of its  singularity  set, in space-time,   is contained in a set whose one-dimensional Hausdorff measure is zero. This means that the singularity set, in space-time, is  ``less" than one-dimensional.

\bu$\,$  Seregin and Sverak \cite{SESV} :

 Smooth solutions of the $3d$ Navier-Stokes equations remain smooth for as long as either the  pressure is bounded from below, or the Bernoulli pressure  $\frac{|v|^2}2 +p$ is bounded from above.

 For instance what \cite{CKN} says is that the velocity of the fluid cannot be too big, say $|u|>M$, except in regions which are small with respect  to $M$.
The result of \cite{SESV} can also be interpreted by saying that for as long as there are no cavities in the flow the solution  can be uniformly controlled. Of course, the minute cavities are formed   the Navier-Stokes equations are  no longer valid as a physical model.


\subsection{ Zero viscosity limit  of solutions of Navier-Stokes in the absence of  boundary}\label{easywlim}

As mentioned in the previous section, the introduction of the viscous term, i.e., the Navier-Stokes equations, is essential in order to  take into account boundary effects. Conversely, in the absence of physical boundaries, one has the following proposition: the behavior, as the viscosity  $\nu\rightarrow 0 $, of Leray-Hopf weak solutions of the Navier-Stokes equations is related to the behavior of   the solutions of the Euler equations according to
\begin{proposition}\label{viscolimi}
1. For given initial data, $u_0 \in L^2(\R^d)$,  Leray-Hopf weak solutions of the Navier-Stokes equations are  bounded (uniformly with respect to $\nu$) in $L^\infty(\R_t; L^2(\R_x^d))\,.$

\noindent 2. Any weak limit, $\overline u$, as $\nu \to 0$, of such a  family of Leray-Hopf weak solutions of Navier-Stokes equations  is a dissipative solution of   Euler equations.

\noindent 3. Assume that, with $u_0$ as initial data,  Euler equations have a smooth solution $ u(t)$, on the time interval  $[0, T]$. Then  $\overline u (t)=u(t)$, for all $t \in [0, T]$.
\end{proposition}\label{prekato}
{  {\begin{remark}
The issue in the above proposition is that, in the absence of boundary effects (with boundary the situation  is completely different, see section \ref{vbdef}, below), any weak limit of solutions to the Navier-Stokes equations is a dissipative solution to the Euler equation. Then the next issue does not concern  Navier-Stokes equations but the dissipative solution of the Euler equations. In the presence of a smooth solution (with the same initial data) of  Euler equations point 2 of Poposition \ref{disssolr} implies that this weak limit is the smooth solution of the Euler equations. Hence, it is unique and nothing spectacular happens as described  with more details in the Remark \ref{conclusion0}, below.
\end{remark}}}
{\bf Proof  of the proposition \ref{prekato}:}
 By definition, a Leray-Hopf weak  solutionsb  $u_\nu(x,t)$  to Navier-Stokes equations satisfies the energy inequality (here without  boundary terms, since there is no physical boundary):
\begin{equation}
\frac12\int_{\R^d} |u_\nu(x,t)|^2dx +\nu\int_0^t\int_{\R^d} |\nabla u_\nu(x,s)|^2dxds \le \frac12\int_{\R^d}  |u_0(x)|^2dx\,.
\end{equation}
Therefore, $u_\nu$ is uniformly bounded in $L^\infty(\R_t; L^2(\R_x^d))$  and, modulo extraction of a subsequence, converges weak$-*$ in this space to a limit $\overline u$.  Proceeding as in section \ref{wsol} (for rigorous justification of all the steps below one can follow the arguments introduced in \cite{Serrin63}), for any divergence-free smooth test function $w(x,t)$ satisfying
\begin{equation}
\del_t w +P_\sigma((w\cdot \nabla) w)=E(w)\,,
\end{equation}
one obtains the estimate:
\begin{equation}
\begin{aligned}
&\int_{\R^d}|u_\nu (x,t)-w(x,t)|^2 dx+2\nu\int_0^t \int_{\R^d} |\nabla u_\nu(x,s)|^2dxds \le \\
&2 \int_0^t\int_{\R^d} |(E(w(x,s))\cdot  (u(x,s)-w(x,s)))|dxds+2\nu\int_0^t \int_{\R^d} (\nabla u_\nu(x,s)\cdot  \nabla w(x,s)) dx ds
\\
&+2\int_0^t \int_{\R^d}|( u_\nu(x,s)-w(x,s) S(w) (u_\nu (x,s)-w(x,s))|dx ds\\
&+ \int_{\R^d} |u_0(x) -w(x,0)|^2dx\,.
\end{aligned}
\end{equation}
Thanks to  Gr\"onwall Lemma one has:
\begin{equation}
\begin{aligned}
&\int_{\R^d} |u_\nu (x,t)-w(x,t)|^2\le  \Bigg( \int_{\R^d} |u_0(x)-w(x,0)|^2dx \Bigg ) e^{2\int_0^t|\!|S(w)|\!|_\infty(s) ds} \\
&+2 \int_0^t \Bigg[\int_{\R^d}|(E(w(x,s))\cdot ( u_\nu(x,s)-w(x,s)))|e^{2 \int_s^t|\!|S(w)|\!|_\infty(\tau) d\tau}dx \Bigg]ds\\
&+2\nu\int_0^t \Bigg[\int_{\R^d} (\nabla u_\nu(x,s)\cdot \nabla w(x,s) ) dx  e^{2 \int_s^t|\!|S(w)|\!|_\infty(\tau) d\tau}\Bigg ]ds\,.\label{dissolns}
\end{aligned}
\end{equation}
Thus, by virtue of the energy estimate:
$$
2 \nu\int_0^T\int_{\R^d} |\nabla u_\nu (x,s)|^2dxds \le  \int_{\R^d}  |u_0(x)|^2dx
$$
there exists a subsequence of $u_\nu$, which will also be labeled $u_\nu$,
such that ${\sqrt \nu \nabla u_\nu} $ converges weakly  in $L^2([0,T]; L^2(\R_x^d))$. As a result one has:
$$
\lim_{\nu \rightarrow 0} 2\nu\int_0^t \Bigg[\int_{\R^d}  (\nabla u_\nu(x,s)\cdot \nabla w(x,s) )dx  \, e^{2 \int_s^t|\!|S(w)|\!|_\infty(\tau) d\tau}\Bigg ]ds=0\,.
$$
Moreover, as noticed above, the weak convergence gives
\begin{equation*}
\begin{aligned}
&\int_{R^d} |\overline u (x,t)-w(x,t)|^2 \le \liminf_{\nu\rightarrow 0} \int_\Omega |u_\nu (x,t)-w(x,t)|^2\\
&\hbox{ and } \\
&\int_0^t \Bigg [\int_{\R^d}|(E(w(x,s))\cdot (\overline u(x, s)-w(x,s)))|\,e^{2 \int_s^t|\!|S(w)|\!|_\infty(\tau) d\tau}dx \Bigg]ds
\\
&=\lim_{\nu\rightarrow 0} \int_0^t\Bigg [\int_{\R^d}|(E(w(x,s))\cdot ( u_\nu(x,s)-w(x,s)))|dx\,e^{2 \int_s^t|\!|S(w)|\!|_\infty(\tau) d\tau}\Bigg ]ds\,.
\end{aligned}
\end{equation*}
Consequently, one has:
\begin{equation}
\begin{aligned}
&\int |\overline u  (x,t)-w(x,t)|^2\le \Bigg( \int |u_0(x) -w(x,0)|^2dx \Bigg)e^{2\int_0^t|\!|S(w)|\!|_\infty(s) ds} \\
&+2 \int_0^t\Bigg [\int_{\R^d}|(E(w(x,s)) \cdot  (\overline u  (x,s)-w(x,s)))dx |e^{2 \int_s^t|\!|S(w)|\!|_\infty(\tau) d\tau} \Bigg ]ds\,.\\
\end{aligned}
\end{equation}
\begin{remark}\label{conclusion0}
1. The most important remark following this proof is the fact that in the presence of a smooth solution of  Euler equations any weak limit will coincide with this smooth solution
and since the smooth solution conserves   energy  one has:
\begin{equation}
\begin{aligned}
&\lim_{\nu\rightarrow 0}\Big( \nu\int_0^t\int_{\R^d} |\nabla u_\nu(x,s)|dxds \Big)\\
& \le   \frac12\int_{\R^d} |u_0(x)|^2dx-
\lim_{\nu\rightarrow 0} \frac12\int_{\R^d} |u_\nu(x,t)|^2dx\\
&\le  \frac12\int_{\R^d} | u_0(x)|^2dx -\frac12\int_{\R^d} |u(x,t)|^2dx=0
 \end{aligned}
 \end{equation}
and therefore there is no anomalous dissipation of energy (the  energy dissipation rate which appears in Kolmogorov theory of turbulence is 0.)

\noindent 2. In fact, as said above, regularity is known to propagate for smooth solution of the $2d$ Euler equations, therefore, in $2d$  with smooth initial data the energy will  always be conserved and there will be no anomalous dissipation of energy (or enstrophy). The appearance of singularities, after some (long?) time, for the $3d$ Euler equations with smooth initial data is still a challenging open problem. Therefore, it may happen (even with smooth initial data) that, after some long time,  anomalous dissipation of energy appears in the limit, as $\nu\rightarrow 0$. Thus,  anomalous dissipation  would be a criterion for blow up of the solutions of  $3d$ Euler equations.

\noindent 3. Moreover, returning to the contributions of De Lellis and Sz\'ekelyhidi, one should keep in mind that there exists in $2d$ or $3d$ an infinite set of initial data, $u_0\in L^2$, such that for each $u_0$  there is an infinite set of corresponding dissipative solutions for  Euler equations. In this situation the above consideration does not apply. For such initial data there may be anomalous dissipation of energy. In this case, since the fact of being dissipative does not imply the uniqueness of the solution, one should consider other types of uniqueness criteria, or selection mechanism for uniqueness. A natural candidate for unique solution of Euler equations would be a viscosity limit solution, namely the limit solution, $\bar{u}$,  of  a family of solutions, $u_\nu$, of  Navier-Stokes equations, when the viscosity $\nu \to 0$. This assertion is not proven with full generality, but some examples built with the shear flow, cf.\cite{BTW}, indicate that this should  be the case.
\end{remark}

\begin{remark}
   Extension of the Proposition \ref{viscolimi}, assuming more regularity on the initial data,
 can be found in   \cite{masmoudi} (see also \cite{CF88}). Assuming that Euler equations have a smooth solution on the time interval  $[0, T]$, it is shown that the same is true for the solution of Navier-Stokes equations with the same initial data on the same time interval for $\nu $ small enough. Moreover, the rates of convergence in higher norms are estimated (see also \cite{Linshiz-Titi} for similar results concerning the Euler$-\alpha$ regularization model.)
 \end{remark}

\subsection{Viscosity limit with boundary effects}\label{vbdef}
 Unlike the previous section the emphasis here is put on  the viscosity limit, in  both the $2d$ and the $3d$ cases, in the presence of physical boundaries. In fact, in as far as our understanding of  this problem goes, there is not much difference between  the $2d$ and the $3d$ cases. Thus,  it will be assumed that the solution of  Euler equations (with given fixed and smooth initial data $u_0$) is smooth for $t\in [0,T]$. In the $2d$ case the global regularity of Euler equations, i.e. for all $T>0$,   was proven by Wolibner in 1933 \cite{wolibner} (and independently at the same time by E.~H\"older \cite{Holder}),  and it will be assumed in the present discussion for the $3d$ case as well.  As said above, the introduction of   Navier-Stokes equations allows to take into account additional boundary effects (in particular the production and shedding of vorticity  near the boundary). One  may link these effects with anomalous dissipation of energy, as $\nu \to 0$, and conclude that turbulent regimes which are observed in the fluid flow,  even far from the boundary, where one could talk of {\it homogenous and isotropic turbulence,}   is not a smooth solution of   Euler equations.

In fact, the number of boundary conditions required for Navier-Stokes equations is  $3$, and  this number is  reduced  to $1$ in the case of Euler equations, i.e., $u\cdot\vec n=0$, generating a boundary layer near the boundary . However, since the problem is nonlinear, this boundary layer may not be confined to the boundary and instead  it may propagate, by means of the nonlinear advection term,  inside the bulk of the flow.

Here also dissipative solutions to the Euler equations will be used; but  an important modification  will be introduced in their definition. In the present case, smooth solutions of Euler equations may have nonzero tangential component at the physical  boundary, thus comparing them  to test functions that have compact support in the domain (i.e. without interaction with the boundary) is not enough to characterize the limit of Navier-Stokes solutions, as the viscosity $\nu \to 0$. The modified  definition of dissipative solutions for the Euler equations is:

\begin{defi}\label{Def-up-to-Boundary}
Consider the Euler equations in a domain $\Omega$, with nonempty  boundary $\del \Omega \neq \emptyset $.  A vector field  $u$ is  called  a dissipative solution up to the boundary, for  the Euler equations,  if for every smooth   divergence-free test vector field, $w$, which is tangent  to the boundary $\del \Omega$ (i.e., $w\cdot \vec n=0$ on  $\del \Omega$),  and which has  bounded support in $\Omega$, one has the estimate:
\begin{equation}
\begin{aligned}
&\frac 12 \int_\Omega |u(x,t)-w(x,t)|^2\le  \int_0^t\int_\Omega |(E(w(x,s)), u(x,s)-w(x,s))|dxds\\
&+\int_0^t \int_\Omega |( u(x,s)-w(x,s) S(w) u(x,s)-w(x,s))|dx ds\\
&+ \frac 12 \int_\Omega |u(x,0)-w(x,0)|^2dx\,, \label{dissolbd}
\end{aligned}
\end{equation}
where $E(w)$ is the residual as in equation (\ref{Res}).
\end{defi}
This is the same criterion as in (\ref{dissol}), but with a larger class of test functions, $w$.
In particular, if there exists a smooth solution $u(x,t)$ with initial data $u_0(x)$, any dissipative solution up to the boundary with the same initial data does coincide with this smooth solution.
Such a result  would not be valid  if only test functions with compact support in $\Omega$ were used  in the criterion (\ref{dissol}).

\begin{theorem}\label{kkato}

 In dimensions $2$ or $3$, for any $\nu$ dependent family of Leray-Hopf weak solutions, on the time interval $[0,T]$, of the Navier-Stokes equations
\begin{equation*}
\del_t u_\nu + u_\nu\cdot \nabla u_\nu -\nu \Delta u_\nu+ \nabla p_\nu=0\,,
\end{equation*}
 in a domain $\Omega$, with boundary conditions
\begin{eqnarray}
&&{u_\nu \cdot \vec n =0\,, \nu ( \del_{\vec n}u_\nu  + C(x) u_\nu )_\tau+\lambda(\nu,x)u_\nu =0\,\hbox{ on } \del \Omega\,, \label{boundcond2}}\\
&& \hbox{and} \,\,\lambda(\nu,x)\ge 0 \,,\,\, C(x)\in C( \R^n; \R^n)\,; \label{coeff2}
\end{eqnarray}
and with fixed initial data $u_0$,  one has the following  assertions:

1. If $u_\nu$ converges weak$-*$, in $L^\infty ((0,T);L^2(\Omega))$, to a function $u(x,t)$, which is, for all $t\in (0,T)$, a smooth solution of the  Euler equations, then  $u_\nu$ converges strongly to $u$, and it does not dissipate energy in the limit as $\nu \rightarrow 0$; namely,
\begin{equation}
\epsilon(\nu,T)= \lim_{\nu\rightarrow 0}\int_0^T\!\!( \nu \int_\Omega |\nabla u_\nu|^2dx +\int_{\del \Omega} \lambda(\nu)|u_\nu(x,t)|^2d\sigma )dt =0\,.\label{enbound}
\end{equation}

2.  If $(\nu \frac{\del u_\nu}{\del \vec n})_\tau \rightarrow 0  $, as $\nu \to 0$, even in a very weak sense (for instance, in the sense of distribution in ${\mathcal D'}(\del\Omega \times(0,T))$), then any weak limit, $\overline u$, of the sequence $u_\nu$ is a dissipative solution up to the boundary of the Euler equations.

3. (Kato criterion) If the sequence $u_\nu$ satisfies the estimate:

\begin{equation}
\lim_{\nu\rightarrow 0} \int_0^T \nu \int_{\{d(x,\del\Omega) < c\nu\}\bigcap \Omega}|\nabla u_\nu(x,t)|^2dxdt =0\,, \label{katocri}
\end{equation}
then the hypothesis of the point 2  above is satisfied, and therefore any weak limit of the sequence $u_\nu$ is a dissipative solution up to the boundary of Euler equations.
\end{theorem}
{\bf Proof}
Below, we will sketch  some of the main ideas of the proof.

Equation (\ref{enbound}) represents the dissipation of energy  over the time $[0,T]$, the term
$\int_0^T\!\!( \nu \int_\Omega |\nabla u_\nu|^2dx)dt $ represents the dissipation of energy due to the fluid viscosity,  and   $\int_0^T(\int_{\del \Omega} \lambda(\nu)|u_\nu(x,t)|^2d\sigma )dt$ represents the dissipation of energy due to  the  friction of the fluid flow with the boundary. In the special case of the Dirichlet boundary condition for the Navier-Stokes equations this second term does not appear  because the velocity field,  $u_\nu$, is zero on the boundary.

Point 1 follows from the energy estimate
\begin{equation*}
\begin{aligned}
\frac 12\int_\Omega |u_\nu(x,t)|^2dx + &\nu\int_0^T \int_\Omega |\nabla u_\nu|^2dx +\int_{\del \Omega} \lambda(x)|u_\nu(x,t)|^2d\sigma\\
&\le \int_\Omega |u_0|^2dx	 +o(1)\,,
\end{aligned}
\end{equation*}
where  $o(1)$ stands for a term that tends  to zero, as $\nu \to 0$.
The above energy inequality is obtainedin a similar way as for the Leray-Hopf weak solutions of the Navier-Stokes system;  however, this time using the   boundary conditions (\ref{boundcond2})-(\ref{coeff2}). For this reason we also call weak solutions of Navier-Stokes equations that satisfy the above energy inequality Leray-Hopf weak solutions.

Observe, thanks to the    weak$-*$ convergence in $L^\infty ((0,T);L^2(\Omega))$ of $u_\nu$ to a smooth solution $u$ of Euler equations,  that one has:
\begin{equation}
\begin{aligned}
\frac 12\int_\Omega |u_0(x))|^2dx=\frac 12\int_\Omega |u(x,t)|^2dx&=\frac 12\int_\Omega |u(x,t)|^2dx\\
&\le \lim_{\nu \rightarrow 0}\frac 12\int_\Omega |u_\nu (x,t)|^2dx\,.
\end{aligned}
\end{equation}

Concerning point 2, the main difference between  the  problem without physical boundary and the current problem with boundary comes from the, only uncontrolled, extra term in the computation of
$$
-\nu\int_0^T\int_{\Omega}( \Delta u_\nu \cdot (u_\nu-w))dxdt\,,
$$
which, thanks to the boundary conditions $w\cdot \vec n=0$ and  (\ref{boundcond2})-(\ref{coeff2}),  gives a boundary term of the form
\begin{equation}
\nu \int_0^T \int_{\del\Omega}((\frac{\del u_\nu}{\del \vec n})_{_\tau }\cdot w_\tau) d\sigma dt\,. \label{bad}
\end{equation}
Thus to establish point 2 one needs to show that the above term goes to zero, as $\nu \to 0$. By virtue of  the density of the set of {  {smooth }}test functions, $w$,  in Definition \ref{Def-up-to-Boundary}  (see  estimate (\ref{dissolbd})), it is enough to consider  in  the criterion stated in point 2, and to test  it for,  $C^\infty$ test functions.

Next, we show point 3.  For any smooth $w(x,t)$ tangent to $\del \Omega$ for all $t\in (0,T]$,  one  constructs
(following Kato \cite{KA}) a family  $w_\nu(s,\tau,t)$ (here $\tau$ and $s$ are  the local tangential and normal coordinates at $\del \Omega$, respectively), for small $\nu >0$,  of smooth divergence-free test functions, which are tangent to the boundary, $\del \Omega$, and which also satisfy  the following properties: For every $\nu >0$, $w_\nu$ coincides with $w$ on  $\del \Omega \times (0,T)$, and its  support is contained in
$$
\Omega_\nu \times (0,T)=\{x\in \Omega:\, d(x,\del\Omega) <\nu\}\times (0,T) \,.
$$
Moreover,   $w_\nu$ satisfied  the  estimates
\begin{eqnarray}
|\del_t w_\nu |_{L^\infty}\le C\,,\quad|\nabla_{\tau}w_\nu |_{L^\infty}\le C\,, \quad |\del_s w_\nu |_{L^\infty} \le \frac{C}{\nu}\,. \label{Approx-w}
\end{eqnarray}
   Next, we show that   (\ref{katocri}) implies   $\lim_{\nu \to 0}(( \nu\frac{\del u_\nu}{\del \vec n}){_{\tau }},w_\tau)_{L^2(\del\Omega\times (0,T))} =0$.   Observer that  $((\nu \frac{\del u_\nu}{\del \vec n}){_{\tau }},w_\tau)_{L^2(\del\Omega\times (0,T))}=(( \nu\frac{\del u_\nu}{\del \vec n}){_{\tau }},(w_\nu)_\tau)_{L^2(\del\Omega\times (0,T))}$, and hence this expression  can be  deduced from Navier-Stokes equations, subject to the boundary conditions (\ref{boundcond2})-(\ref{coeff2}), according to the formula:
\begin{equation}
\begin{aligned}
&0=\int_0^T(0,w_\nu) dt=\int_0^T((\del_t u_\nu +\nabla\cdot ( u_\nu \otimes u_\nu)-\nu \Delta u_\nu +\nabla p_\nu), w_\nu)dt = \\ & -\int_0^T(u_\nu,\del_t w_\nu)dt
 - \int_0^T(( u_\nu \cdot \nabla) w_\nu,  u_\nu) dt \\ & +\nu\int_0^T(\nabla u_\nu,\nabla w_\nu)dt  -(\nu \del_{\vec n } u_\nu,w_\nu)_{L^2(\del\Omega\times (0,T))}\,,
\end{aligned}
\end{equation}
which implies the estimate:
\begin{equation}\label{Est1}
\Bigg|(\nu \del_{\vec n } u_\nu, w_\nu)_{L^2(\del\Omega\times (0,T))}\Bigg|=\Bigg|\int_0^T(( u_\nu \cdot \nabla) w_\nu,  u_\nu) dt\Bigg|+o(1)\,.
\end{equation}
Eventually, using (\ref{Approx-w}), the control of the energy and the Poincar\'e inequality, one obtains
\begin{equation}\label{Est2}
\Bigg|\int_0^T(( u_\nu \cdot \nabla) w_\nu,  u_\nu) dt\Bigg| \le C \int_0^T\int_{\Omega_\nu} \nu |\nabla u_\nu|^2dxdt\,.
\end{equation}
Therefore, thanks to (\ref{katocri}), (\ref{Est1}) and (\ref{Est2}) one concludes
$$
\lim_{\nu \to 0}(( \nu\frac{\del u_\nu}{\del \vec n}){_{\tau }},w_\tau)_{L^2(\del\Omega\times (0,T))} =0\,.
$$
This  completes  the proof of  point 3.
\begin{remark}

The boundary conditions (\ref{boundcond2})-(\ref{coeff2})
together with the hypothesis $\lim_{\nu \rightarrow 0}\lambda(\nu)=0$ imply that any weak limit of the family $u_\nu$, of Leray-Hopf weak solutions of the Navier-Stokes equations, is a dissipative solution up to the boundary of Euler equations.
{  {The hypothesis of  point 2 involves only $\nu$ times the tangential component of the normal derivative of the velocity. This point of view can be combined with the hypothesis of  point 3, and this was done by Wang \cite{wang}. The condition proposed by Wang (which at the end of the day will be equivalent to all the others) involves only   the tangential component of the normal derivative of the velocity, but in a layer slightly thicker  than the Kato layer, which is the well-known thickness of viscous sublayer as predicted by  von Karman (1930) logarithmic profile.
In situation when the convergence to the regular solution of  Euler equations is insured it is possible to give many details of the rate of this convergence, the nature of the boundary layer and so on. There is a huge literature along these lines, to cite few see, for example,  \cite{bdv1}, \cite{da-Veiga-Crispo}, \cite{zhouping} and references therein. }}
\end{remark}
\begin{remark}
When $\lambda(\nu)$ does not go to zero, in particular,  in the most extreme case (the Dirichlet boundary conditions) when   $\lambda(\nu)=\infty$,  there is no information on the nature of the limit, at this level of generality.

However, one can underline the interpretation of the above theorem by stating that any of the   following facts:

\bu There is no  production of  vorticity at the boundary.

\bu There is no anomalous dissipation of energy in  a layer of size $\nu$, at the boundary, and, of course, no anomalous dissipation of energy in the whole fluid flow domain.

\bu There  is conservation of energy as  $\nu\rightarrow 0\,,$

\noindent
implies that the zero viscosity limit   does not exhibit  any turbulent behavior. This is meant in the sense  that  the limit  solution of Navier-Stokes equations, as $\nu \to 0$,  is  as smooth as the corresponding solution of Euler equations,  and that no nontrivial Reynolds stresses tensor (in the sense of coarse graining) appears at the limit.

On the other hand, if one of the above statement is not true,  due to boundary effects,  one cannot say that the limit  is a dissipative solution up to the boundary  of Euler equations in contrast with the case without boundaries.

Nontrivial Reynolds stresses tensor (in the sense of coarse graining)  may occur  either when the weak limit, $\overline u $, as $\nu \to 0$,  is not a weak solution of   Euler equations, or when$\overline u $ is a weak solution of Euler equations that  is obtained as the limit  of Navier-Stokes equations with energy decay, i.e., $\overline u $  is a weak admissible solution of Euler equations. However, even so, some concentration of oscillations (turbulence) near the boundary may prevent it to be a dissipative solution up to the boundary of Euler equations. This pathologic type of weak solutions of Euler equations can in fact be constructed following \cite{BS}.
\end{remark}

\subsection{More about boundary effects and Prandtl equations}
{  {In this section we follow very closely our presentation (cf. \cite{Bardos-Titi-survey}) concerning the Prandtl equations.
What essentially prevents the existence of more definite results, in the general case, is that what would be a boundary layer, in the zero viscosity limit, may not remain near the boundary and will actually detach from the boundary due to its rolling up into vortices which will, by the advection term, penetrate the bulk of the fluid flow forming the wake of the solid body.}}

However, this fact does not prevent more specific theoretical or practical studies. If one assumes, as a hypothesis, that the variation of the tangential component, and  the size of the normal component, of the velocity of the fluid, near the boundary, are of the order of $\epsilon^{-1}$ and $\epsilon$, respectively,  (with $\epsilon=\sqrt \nu$), then by a change of scale near this boundary one may use the Prandtl
 equations (cf. \cite{PRA} and for more recent discussion concerning the relation between the issue of $0$ viscosity limit and the Prandlt equations \cite{wang2}) We recall below the planar  Prandtl equations for the case of a flat boundary:
start from the  Navier-Stokes equations in the half-plane,  $\{(x_1,x_2)\in \R^2:\,x_2>0\}$, with the no slip boundary condition, $u_\nu(x_1,0,t)\equiv 0$,
\begin{eqnarray}
&&\del_t u_1^\nu -  \nu\Delta u_1^\nu + u_1^\nu\del_{x_1} u_1^\nu
+u_2^\nu\del_{x_2} u_1^\nu +\del_{x_1} p^\nu=0\,,
\\
&&\del_t u_2^\nu - \nu \Delta u_2^\nu + u_1^\nu\del_{x_1} u_2 ^\nu+
u_2^\nu\del_{x_2} u_2^\nu
+\del_{x_2} p^\nu=0\,,
\\
&&\del_{x_1} u_1^\nu+\del_{x_2} u_2^\nu =0\,,
\\
&&
u_1^\nu(x_1,0)=u_2^\nu(x_1,0)=0  \hbox { for }  x_1\in \R\,,
\end{eqnarray}
then introduce the scale  $\epsilon=\sqrt \nu$ and the change of variable which corresponds to a boundary layer in a parabolic problem, taking into  account that the normal component of the velocity on the boundary, $u_2$, remains $0$:
\begin{eqnarray}
&&X_1=x_1, X_2=\frac{x_2}{\eps},    \\
&&
\uT_1(x_1,X_2)=u_1(x_1,x_2),
\uT_2(x_1,X_2)=\eps u_2(x_1,x_2)\,.\end{eqnarray}
Going back to the original notation  $(x_1,x_2)$ and letting $\epsilon$ go to zero, one obtains formally the system:
\begin{eqnarray}
&&  \uT_1(x_1,0,t )-U_1(x_1,t)=0\,,  \label{Eulerext}\\
&& \del_{x_2}\tilde p(x_1,x_2)=0\Rightarrow  \tilde p  (x_1,x_2,t) =\tilde P (x_1,t)\label{Eulerext2}\\
&&\del_t \uT_1  - \del_{x_2}^2 \uT_1 + \uT_1 \del_{x_1} \uT_1
+\uT_2 \del_{x_2} \uT_1   =-\del_{x_1} \tilde P(x_1,t)\,, \label{Eulerext3}\\
&&
\del_{x_1} \uT_1 +\del_{x_2} \uT_2  =0\,,  \\
&&
\uT_1(x_1,0,t)=\uT_2(x_1,0,t) =0 \hbox { for } x_1\in \R\,,
\\
&&\lim _{x_2\rightarrow \infty} \uT_1(x_1,x_2)=0\,.
\end{eqnarray}

Equations  (\ref{Eulerext}) 	and  (\ref{Eulerext2}) describe the fact that the boundary layer solution will match a solution of   incompressible Euler equations.  $U_1(x_1,t)$
and $\tilde P(x_1,t)$ are the tangential velocity and the pressure, respectively,  on the boundary, $x_2=0$,  of the solution to the Euler equations with the same initial data.
Finally, the solution is written as a sum of   regular functions: the boundary layer term, as  computed above, and the solution of Euler equations, denoted $u_{\mathrm {int}}(x_1,x_2,t)$
\begin{equation}
 \Bigg(\begin{matrix}{}u_1^\nu(x_1,x_2,t) \cr
u_2^\nu(x_1,x_2,t)
\end{matrix}  \Bigg)
=\Bigg(\begin{matrix}{} \uT_1^\nu(x_1,\frac{x_2}{\sqrt{\nu}},t) \cr
\sqrt{\nu} \uT_2^\nu(x_1,\frac{x_2}{\sqrt{\nu}},t)\end{matrix}\Bigg)+ u_{\mathrm {int}}(x_1,x_2,t)\label{couchlim}\,.
\end{equation}

\begin{remark}
 The validity of the Prandtl approximation is consistent with the Kato's criterion  (point 3 of Theorem  \ref{kkato}) because with (\ref{couchlim}) one has:
\begin{equation}
\nu\int_0^T\!\!\!\int_{\Omega\cap \{d(x,\del \Omega)\le c \nu\}} |\nabla \wedge u_\nu(x,s)|^2dxds\le C\sqrt{ \nu}\,.
\end{equation}
\end{remark}
\begin{remark}
At least one example showing that the Prandtl expansion cannot be always valid has been constructed by Grenier \cite{GR}. Solutions are considered in the domain
$${\R_{x_1}}/{\mathbf Z}\times\R_{x_2}^+\,.$$
Grenier starts with a solution $u^\nu_{ref}$ of the pressure-less Navier-Stokes equations given by:
\begin{eqnarray}
&&u^\nu_{ref}=(u_{ref}(t,\frac{y}{\sqrt{\nu}}),0)\,,\\
&&\del_t u_{ref} -\del^2_{Y}u_{ref}=0\,,
\end{eqnarray}
where $Y=\frac{y}{\sqrt{\nu}}$\,.
Using a convenient and explicit  choice of the function $u_{ref}$,  with some sharp results on   instabilities, a solution of  Euler equations of the form
\begin{equation}
\tilde u =u_{ref} +\delta v + O(\delta^2 e^{2\lambda t} ) \quad {\rm  for }\quad  0<t<\frac 1{\lambda |\log \delta|}
\end{equation}
is constructed. It is then shown that the vorticity generated at the boundary, by the solution of  Navier-Stokes equations (with the same initial data), is too strong to allow the convergence of the Prandtl expansion.
One should observe that, once again, this is an example with infinite energy. It would be interesting to see if such example could be modified to enter the class of finite energy solutions and analyze how it would violate the Kato's criterion (\ref{katocri}).
\end{remark}
It is  important to observe that the mathematical properties of Prandtl equations exhibit the pathology of a situation based on the hypothesis that the boundary layer does not detach.
First, one can prove the following
\begin{prop} Let $T>0$ a finite positive time, and let $(U(x,t),P(x,t))\in C^{2+\alpha}(\R_t\times (\R_{x_1}\times \R_{x_2}^+))$ be a smooth solution of the $2d$ Euler equation satisfying for $t=0$ the compatibility condition $U_1(x_1,0,t)=U_2(x_1,0,t)=0$ (only the boundary condition $U_2(x_1,0,t)=0$ is preserved by the Euler's dynamics). Then  the following facts are equivalent:

1. with initial data $\uT=0\,,$ the boundary condition $\uT_1(x_1,0,t)=U_1(x_1,0,t)$ in (\ref{Eulerext}) and the  right-hand side in (\ref{Eulerext2})  given by
$\tilde P(x_1,t)= P(x_1,0,t)\,,$  Prandtl equations have a smooth solution for $0<t<T\,.$

2. The solution $u_\nu(x,t)$,  with initial data $u_\nu(x,0)=U(x,0)$ and with no-slip boundary condition, at the boundary $x_1=0$, converges in $C^{2+\alpha}$ to the solution of  Euler equations, as $\nu \to 0$.
\end{prop}
The fact that the hypotheses 1.  (and 2.), above,  may be violated for some $t$ is related to the appearance of a detachment zone and the generation of vorticity at the boundary (which may eventually generate turbulence in the bulk of the fluid flow). This is well illustrated in the analysis of   Prandtl equations   written in the following simplified form:
\begin{eqnarray}
&&\del_t \uT_1  - \del_{x_2}^2 \uT_1 + \uT_1 \del_{x_1} \uT_1
+\uT_2 \del_{x_2} \uT_1   =0\,, \label{pes1}
\del_{x_1} \uT_1 +\del_{x_2} \uT_2  =0Ê\,,
\\
&&\uT_1(x_1,0)=\uT_2(x_1,0) =0 \hbox { for } x_1\in \R\,,\label{pes2}
\\
&&\lim _{x_2\rightarrow \infty} \uT_1(x_1,x_2)=0\,, \\
&& \tilde u_1(x_1,x_2,0)=\tilde u_0(x_1,x_2)\,.
\end{eqnarray}
Regularity in the absence of detachment corresponds to a theorem of  Oleinik \cite{OL}. She proved that global smooth solution does exist  if  the initial profile is monotonic, i.e.
\begin{equation}
\uT(x,0)=(\uT_1(x_1,x_2),0)\,, \del_{x_2} \uT_1(x_1,x_2) \not=0\,.
\end{equation}
On  the other hand, initial conditions with ``recirculation properties" leading to a finite time blow up have been constructed by E and Enquist \cite{EE}.  An interesting aspect of these examples is that the blow up generally does not occur  exactly on the boundary, but inside the domain supporting the hypothesis that this blow up is due to the roll up of the boundary layer which increases the velocity gradients.

The above pathology appears in the fact that the Prandtl equations are  highly unstable. This comes from the determination $\uT_2$  in term of $\uT_1$ by the equation:
\begin{equation}
\del_{x_1} \uT_1 +\del_{x_2} \uT_2  =0\,.
\end{equation}
With this observation more recent and systematic studies of these instabilities can be found in \cite{GVD} and \cite{GVTN}.

Therefore it is only with analytic initial data (in fact analyticity  with respect to the tangential variable is enough) that one can obtain (using an   abstract version of the Cauchy-Kowalewsky theorem) the existence of a smooth solution of   Prandtl equations for a finite time and the convergence to the solution of  Euler equation during the same time  (cf.~Asano, Caflisch-Sammartino and Cannone-Lombardo-Sammartino \cite{AS}, \cite{SC2}, \cite{CLS}, {  {and  more recently with less stringent hypothesis on the $x_2$ asymptotic behavior  in \cite{kukavica-vicol}}}.)

\section{Conclusion}

This contribution is organized around the issues of the Cauchy problem, and the relevant boundary effects, in equations of motion of incompressible fluids. It presents some mathematical tools and results (theorems) that have been developed over the last, and the present, century for investigating  certain basic properties of  solutions of   equations of fluid flows. There are many different mathematical approaches, that complement each other,  to address these challenging problems; however, we chose  to focus, in this contribution,    only  on certain mathematical aspects that are concerned with:

\bu$\,$  The relation between the different equations of fluid motion in terms of well identified scaling parameters.

\bu$\,$  The existing results (far from being complete at the time of this writing) about classical solutions of these equations.

\bu$\,$  The notion of weak solutions (as used by Leray in his thesis or in  the sense of distributions, as introduced by Gelfand and Schwartz). In  particular, we recall  that  the notion of weak solutions   (even with conservation or decay of energy) is not precise enough to characterize the convenient (physically relevant)  solutions of the Cauchy problem for Euler equations. This is in spite the fact that the relation between the regularity of solutions of Euler equations  and the non conservation of energy is related to the Onsager conjecture in turbulence theory.

\bu$\,$  The analysis of the zero viscosity limit for  solutions of Navier-Stokes equations, with or without boundary effects. In the absence of boundary effects this convergence is a direct consequence of the existence (or nonexistence)  of smooth solutions, with the same initial data, for the Euler equations; while in the presence of boundary effects the situation is completely different and the non convergence  (following a theorem of T.~Kato) is  related to the anomalous dissipation of energy, in the limit as  $\nu \rightarrow 0$. This anomalous dissipation of energy corresponds, in the Kolmogorov theory of turbulence,  to the non- vanishing limit of the positive mean rate of dissipation of  energy, $\epsilon$, as   $\nu \rightarrow 0$.

\bu$\,$  Lastly,  the Prandtl equations are shortly discussed. We recall that these equations are in general ill-posed and that in the case when they are well-posed, which is the ``laminar regime", the Kato criterion (of  the above mentioned  Kato theorem) is obviously satisfied, and as a result  there is no anomalous dissipation of energy. Beside their importance (in situations  where they are valid) for engineering sciences applications,  the Prandtl  equations became basic models for a systematic mathematical analysis of boundary layer in many other situations.
\vskip0.1cm

{\bf Acknowledgements.}\hspace*{0.03in} {\it We wish to thank Marie Farge and Kai Schneider for the invitation to  the conference ``Fundamental Problems of Turbulence, 50 Years after the Marseille Conference 1961" from which this contribution has originated. Many thanks also to Marie Farge for her careful reading and hence scientific contribution to the present manuscript. The work of E.S.T.~is  supported in part by the   NSF grants DMS-1009950, DMS-1109640,  and DMS-1109645. E.S.T. ~also acknowledges the support of the  Minerva Stiftung/Foundation.}

   \end{document}